\documentclass[11pt]{article}
\usepackage[a4paper,margin=24mm]{geometry}
\usepackage{lmodern,amsmath,amssymb,amsthm,mathtools,booktabs,array,enumitem,longtable}
\usepackage[expansion=false]{microtype}
\usepackage{xcolor}
\usepackage{tikz}
\usetikzlibrary{arrows.meta,shapes.geometric}
\usepackage[font=small,labelfont=bf,labelsep=colon]{caption}
\definecolor{rayteal}{HTML}{1F6E75}
\definecolor{treenavy}{HTML}{22313F}
\definecolor{extraorange}{HTML}{B4522A}
\definecolor{latgray}{HTML}{D4D4D4}
\definecolor{shadeA}{HTML}{E8F1F1}
\usepackage[colorlinks=true,linkcolor=blue!42!black,citecolor=blue!42!black,urlcolor=blue!42!black,bookmarksnumbered=true]{hyperref}
\hypersetup{pdftitle={Random walk on an induced component of the uniform spanning forest},pdfauthor={Xinyi Li, Runsheng Liu and Yuval Peres},pdfsubject={Recurrence and return probabilities}}
\numberwithin{equation}{section}
\newtheorem{theorem}{Theorem}[section]
\newtheorem{proposition}[theorem]{Proposition}
\newtheorem{lemma}[theorem]{Lemma}
\newtheorem{corollary}[theorem]{Corollary}
\theoremstyle{remark}
\newtheorem{remark}[theorem]{Remark}
\newcommand{\Z}{\mathbb Z}
\newcommand{\E}{\mathbb E}
\newcommand{\PP}{\mathbb P}

\newcommand{\res}{R_{\mathrm{eff}}}
\newcommand{\diam}{\operatorname{diam}}
\newcommand{\Past}{\operatorname{Past}}
\newcommand{\Ray}{\operatorname{Ray}}
\newcommand{\WUSF}{\operatorname{WUSF}}
\newcommand{\UST}{\operatorname{UST}}
\newcommand{\LE}{\operatorname{LE}}

\newcommand{\ind}{\mathbf 1}
\newcommand{\1}{\mathbf 1}
\newcommand{\lp}{\ell}
\newcommand{\cR}{\mathcal R}

\title{Random walk on an induced component\\of the uniform spanning forest}
\author{Xinyi Li\thanks{Beijing International Center for Mathematical Research, Peking University, \href{mailto:xinyili@bicmr.pku.edu.cn}{xinyili@bicmr.pku.edu.cn}. } \and Runsheng Liu\thanks{School of Mathematical Sciences, Peking University, \href{mailto:liurunsheng@pku.edu.cn}{liurunsheng@pku.edu.cn}.} \and Yuval Peres\thanks{Beijing Institute of Mathematical Sciences and Applications (BIMSA), \href{mailto:yperes@bimsa.cn}{yperes@bimsa.cn}.}}\date{}
\begin{document}
\maketitle

\begin{center}
    \itshape
    Dedicated to Susu, the newborn son of the first author
\end{center}

\begin{abstract}
Let $T$ be the component of the origin in the uniform spanning forest of $\mathbb{Z}^d$, and consider the induced graph $H$ obtained by retaining every lattice edge with both endpoints in $T$. Lyons, Peres and Sun (2020) showed that $H$ is recurrent if $d \ge 8$ and asked for which $d \in \{5,6,7\}$ this holds.  We prove that simple random walk on $H$ is transient for $d=5$ and recurrent for every $d\ge6$.
We also show that the quenched return probabilities to the origin satisfy $p^H_{2n}(o,o)=n^{-\alpha_d+o(1)}$ a.s., where $\alpha_5=4/3$ and $\alpha_6=1$, while $\alpha_7=4/5$ and $\alpha_d=2/3$ for all $d\ge8$.
\end{abstract}

\section{Introduction}
\subsection{Background and related work}
The uniform spanning tree (UST) of a finite connected graph is chosen uniformly among its spanning trees. Pemantle \cite{Pemantle} constructed its infinite-volume limit on $\Z^d$ and proved that free and wired boundary conditions give the same law, called the uniform spanning forest (USF). The forest is connected for $d\le4$ and has infinitely many components for $d\ge5$. Benjamini, Lyons, Peres and Schramm \cite[Theorem~10.1]{BLPS} completed the description of its ends, proving in particular that every component in $d\ge5$ has a unique infinite simple path from each vertex. Wilson's algorithm \cite{Wilson}, extended to infinite graphs in \cite[Theorem~5.1]{BLPS}, constructs the forest from loop-erased random walks and provides a basic tool for studying its geometry.

Subsequent work quantified this geometry in both the lattice and tree metrics. Benjamini, Kesten, Peres and Schramm \cite{BKPS} developed multipoint connectivity estimates and showed that, in dimensions five through eight, every two distinct components have infinitely many pairs of neighboring vertices. Barlow and J\'arai \cite[Theorem~1.1]{BJ} obtained upper and lower tail bounds for intrinsic ball volumes, establishing volume growth exponent two in every $d\ge5$. Hutchcroft \cite{H} proved sharp size and radius estimates for the past of a vertex, meaning the vertex together with the finite branches separated from infinity when that vertex is removed. These estimates are among the geometric inputs used throughout this paper.

Random walk restricted to the tree edges has also been studied. Morris \cite{Morris} proved that every component of the wired uniform spanning forest on a locally finite graph is recurrent when all edges have unit conductance. In dimensions $d\ge5$, Hutchcroft \cite[Theorem~1.3]{H} proved that every component has spectral dimension $4/3$ and intrinsic walk dimension $3$ almost surely. In particular, the return probability at time $2n$ is $n^{-2/3+o(1)}$ almost surely in every such dimension. Thus the walk on a single tree has the same return exponent throughout the high-dimensional regime.

The lattice embedding suggests a different walk: retain every lattice edge whose endpoints belong to the same tree, and run simple random walk on the resulting induced graph. Its vertex set is unchanged, but the additional edges create cycles and shortcuts, so the tree results alone do not determine its recurrence or heat kernel. Lyons, Peres and Sun initiated the study of these induced graphs \cite{LPS}. They proved recurrence for $d\ge8$ and, for $d\ge9$, linear growth of the effective resistance from the root to the $n$th vertex on the tree ray \cite[Theorems~1.1 and~1.3]{LPS}. They also proved that, conditional on its induced graph, the original tree has the wired uniform spanning-forest law on that graph \cite[Theorem~1.4]{LPS}. Their work left the recurrence question open in dimensions five, six and seven. We resolve these remaining cases and determine how the additional edges affect return probabilities. Lyons, Peres and Sun also asked about resistance to the $n$th ray vertex for $5\le d\le8$. We obtain the resistance scales for exit from a lattice box in dimensions six through eight.

\subsection{Model and main results}\label{sec:model}\label{sec:results}

\begin{figure}[ht]
\centering
\begin{tikzpicture}[scale=0.6,
  occ/.style={circle,fill=treenavy,inner sep=1.5pt},
  emp/.style={circle,fill=latgray,inner sep=1.1pt},
  ray/.style={rayteal,line width=1.1pt},
  br/.style={treenavy,line width=1.1pt},
  ex/.style={extraorange,line width=1.1pt,dashed}]
 \foreach \dx/\ttl/\extra in {0/{The tree $T$}/0,10/{The induced graph $H$}/1}{
  \begin{scope}[xshift=\dx cm]
   \node at (3,4.4) {\ttl};
   \foreach \p in {(0,0),(0,2),(0,3),(1,3),(2,1),(4,0),(5,0),(5,3),(6,0),(6,1),(6,3)} \node[emp] at \p {};
   \draw[ray] (0,1)--(1,1)--(1,2)--(2,2)--(2,3)--(4,3)--(4,2)--(6,2);
   \draw[ray,-{Stealth[length=5pt]}] (6,2)--(6.8,2);
   \draw[br] (1,1)--(1,0)--(3,0)--(3,2);
   \draw[br] (4,1)--(5,1)--(5,2);
   \ifnum\extra=1
    \draw[ex] (2,2)--(4,2); \draw[ex] (3,3)--(3,2); \draw[ex] (4,2)--(4,1); \draw[ex] (3,1)--(4,1);
   \fi
   \foreach \p in {(1,0),(1,1),(1,2),(2,0),(2,2),(2,3),(3,0),(3,1),(3,2),(3,3),(4,1),(4,2),(4,3),(5,1),(5,2),(6,2)} \node[occ] at \p {};
   \node[circle,fill=extraorange,inner sep=2.6pt] at (0,1) {};
   \node[left=3pt] at (0,1) {$o$};
  \end{scope}}
\end{tikzpicture}
\caption{The same vertex set in the tree and induced graph. The teal path follows the ray; dashed edges give the walk additional steps.}\label{fig:slice}
\end{figure}
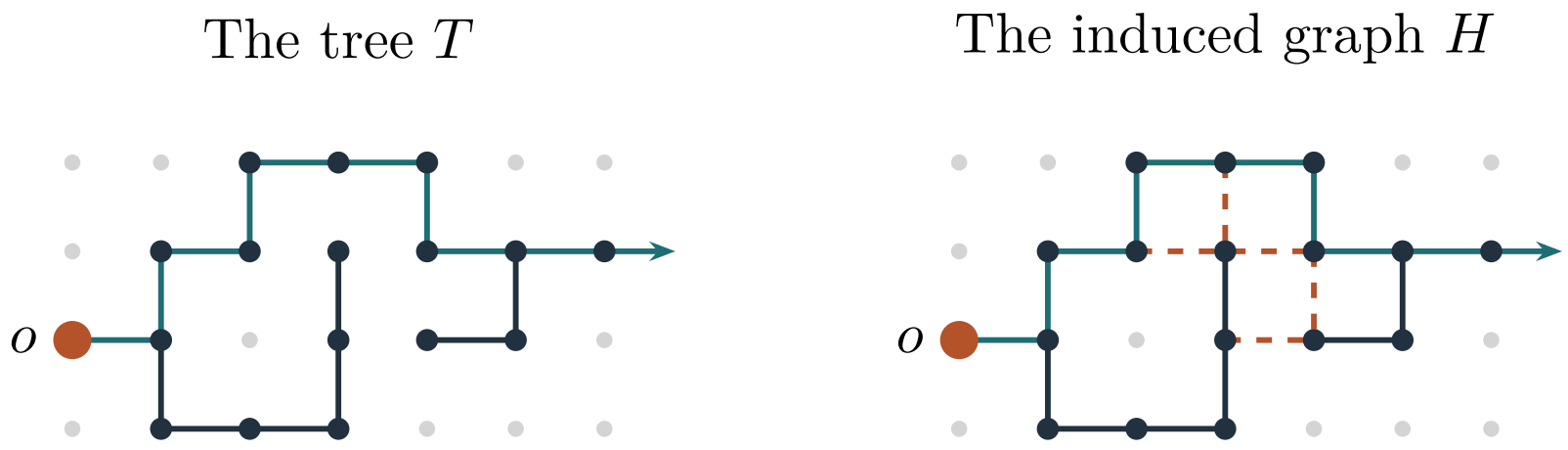

Let $F$ be the wired uniform spanning forest of the nearest-neighbor lattice $G=\Z^d$, $d>4$, and let $T$ be its tree containing $o=0$. We write
\[
 H=G[V(T)],\qquad E(H)=\{\{x,y\}:x,y\in V(T),\ |x-y|_1=1\}.
\]
Here $V(J),E(J)$ denote the vertex and edge sets of a graph $J$, and $|\cdot|_1$ is the $\ell^1$ norm. For each $x\in\Z^d$, similarly write $T_x$ for its forest component and $H_x=G[V(T_x)]$ for the induced graph; thus $T=T_o$ and $H=H_o$. We identify graphs with their vertex sets in set expressions. Every edge has conductance one. Write $\PP,\E$ for the forest law and expectation, and $P_x^H,E_x^H$ for the law and expectation of simple random walk $X=(X_n)_{n\ge0}$ on the fixed graph $H$, started at $x$. Its transition operator is $P_Hf(x)=E_x^H[f(X_1)]$, and $\pi(x)=\deg_{H_x}(x)$ is the reversible measure on each induced component. The one-ended tree $T$ has a unique infinite ray $\gamma=(r_0=o,r_1,\ldots)$. Removing a ray edge separates a finite part of $T$ from infinity; additional edges of $H$ may cross this cut (Figure~\ref{fig:slice}).

Write $\res^J(A,B)$ for effective resistance in a graph $J$ between disjoint vertex sets $A,B$, each wired to one vertex. Singleton sets are denoted by their vertices; for an edge $e=\{x,y\}$, write $\res^J(e)=\res^J(x,y)$. Resistance to infinity is defined along finite wired exhaustions.

Our first theorem locates the recurrence transition between dimensions five and six.
\begin{theorem}\label{thm:main}
In dimension five, almost surely every induced component is transient, and $\E\res^H(o,\infty)<\infty$. For every $d\ge6$, almost surely every induced component is recurrent.
\end{theorem}

For $d\le4$, the forest is connected \cite{Pemantle,BLPS}, so its induced graph is $\Z^d$. Thus the walk on an induced component is transient exactly in dimensions three, four and five.

Beyond recurrence, the additional edges change the return exponent in dimensions five, six and seven, while in dimension eight the fixed-confidence and annealed return estimates have a logarithmic correction. We next state the quantitative results, distinguishing the walk on a fixed induced graph from its average over the forest.

Put $Q_R(x)=x+[-R,R]^d\cap\Z^d$ and $Q_R=Q_R(o)$. For a graph $J$, let $d_J$ be its graph distance, infinite between distinct components, and let $B_J(x,R)=\{y\in V(J):d_J(x,y)\le R\}$ be its closed ball. Write $p_n=p_n^H(o,o)$ with $p_n^H(x,y)=P_x^H(X_n=y)$. Odd returns vanish by bipartiteness. Put
\[
 \cR_R=\res^H(o,H\setminus Q_R),\qquad \cR_\infty=\res^H(o,\infty),
\]
where in $\cR_R$ all vertices of $H$ outside $Q_R$ are wired to one vertex. For nonnegative random $Y_t$ and deterministic $a_t>0$, write $Y_t\asymp_{\PP}a_t$ if, for each $\eta>0$, deterministic $0<c_\eta<C_\eta<\infty$ satisfy $\PP(c_\eta a_t\le Y_t\le C_\eta a_t)\ge1-\eta$ at every sufficiently large deterministic $t$. The symbol $\asymp$ compares deterministic quantities up to positive constants, and $\lp_t=\log(2+t)$ for $t\ge0$.

The next theorem separates fixed-confidence and averaged estimates. In dimensions five and six, we do not determine the order of the mean return probability.
\begin{theorem}[Return probabilities]\label{thm:returns}\label{thm:summary}
The second column gives the fixed-confidence scale of $p_{2n}$ and the third the annealed return $\E p_{2n}$.
\[
\begin{array}{c@{\qquad}c@{\qquad}l}
\toprule
d&p_{2n}\asymp_{\PP}&\text{Annealed return }\E p_{2n}\\\midrule
5&n^{-4/3}&\ge cn^{-4/3}\\
6&n^{-1}&c/n\le\E p_{2n}\le C\lp_n/n\\
7&n^{-4/5}&\asymp n^{-4/5}\\
8&(n\lp_n)^{-2/3}&\asymp(n\lp_n)^{-2/3}\\
\ge9&n^{-2/3}&\asymp n^{-2/3}\\\bottomrule
\end{array}
\]
Almost surely $p_{2n}=n^{-\alpha_d+o(1)}$, where $\alpha_d=4/3,1,4/5,2/3,2/3$ in the five rows, so the quenched spectral dimension $-2\lim_{n\to\infty}\log p_{2n}/\log n$ equals $2\alpha_d$.
\end{theorem}

\emph{Moments and fluctuations.} In dimension five, $\E[p_{2n}^q]\asymp n^{-4q/3}$ for $0<q<3/4$ (Lemma~\ref{lem:momentsfive}). In dimension six, Proposition~\ref{prop:d6-return-bounds} gives $\E[p_{2n}^q]\asymp n^{-q}$ for $0<q<1$ and $\E\sum_{k\le n}p_{2k}\asymp\lp_n$. For every $\varepsilon>0$, eventually almost surely,
\[
 \frac{c_\varepsilon}{n(\log n)^{5+\varepsilon}}\le p_{2n}\le\frac{C_\varepsilon\log n(\log\log n)^{1+\varepsilon}}n.
\]
In dimension eight, $(n\lp_n)^{2/3}p_{2n}$ has uniformly bounded positive moments of every order below $3/2$ (Proposition~\ref{prop:d8-sharp}). In dimensions five, six and seven, $\limsup_n n^{\alpha_d}p_{2n}=\infty$ almost surely, so no deterministic normalization captures $p_{2n}$ up to random constants (Proposition~\ref{prop:return-fluctuations} and Corollary~\ref{cor:return-fluctuation-consequences}).

\emph{Box resistance.} By Lemma~\ref{lem:terminal} and Proposition~\ref{prop:d8-sharp}, the resistance to the exterior of $Q_R$ satisfies
\[
 \cR_R\asymp_{\PP} b_d(R),\qquad \E\cR_R\asymp b_d(R),\qquad
 b_d(R)=
 \begin{cases}
 \lp_R & d=6,\\
 R & d=7,\\
 R^2/\lp_R & d=8.
 \end{cases}
\]

Among dimensions $d\ge5$, the induced graph and the tree have different spectral dimensions exactly for $d=5,6,7$. In dimension eight the spectral dimension is still $4/3$, although the fixed-confidence and annealed return scales in Theorem~\ref{thm:returns} contain a logarithm.

\subsection{How the forest geometry enters}\label{sec:bkps}
The probabilistic input about the forest comes from the work of Benjamini, Kesten, Peres and Schramm \cite{BKPS}, supplemented by the volume bounds of Barlow and J\'arai \cite{BJ} and Hutchcroft's tails for the past of a vertex \cite{H}. Section~\ref{sec:prelim} states the versions used. Three consequences of \cite{BKPS} organize the whole paper, and each is turned into an estimate below.

\emph{Trees are four-dimensional.} Benjamini, Kesten, Peres and Schramm showed that the forest has stochastic dimension four: $\PP(x\leftrightarrow y)\asymp(1+|x-y|)^{4-d}$, together with a multipoint connection bound \cite[Theorem~4.3]{BKPS}. Summing it over $k$-tuples gives $\E|T\cap Q_R|^k\le CR^{4k}$, and Barlow and J\'arai showed that $|T\cap Q_R|$ is rarely much smaller than $R^4$ (Lemma~\ref{lem:forest-volumes}). Two consequences are used throughout. First, the walk on $H$ confined to $Q_R$ has about $R^4$ vertices available; this produces the volume terms in every return bound. Second, dividing the volume $R^d$ by $R^4$ suggests an upper bound of order $R^{d-4}$ on the expected number of trees meeting $Q_R$. Lemma~\ref{lem:packing} proves this packing bound. Using the Barlow--J\'arai volume bound in a box with wired boundary, it also shows that the bound survives conditioning on distant edges, which arises below.

\emph{Neighboring trees and the threshold at dimension six.} Deleting the ray edge chosen by the resistance formula replaces the forest by the forest wired at the origin (Lemma~\ref{prop:deletion}). If $A_o$ is the finite origin tree in the wired-at-$o$ forest, then the mean escape resistance $\E\cR_\infty$ in the original forest equals $r_G\E^{o,\mathrm w}D_o$, where $r_G$ is the lattice escape resistance and $D_o$ counts the trees adjacent to $A_o$. The radius of $A_o$ exceeds $s$ with probability at most $Cs^{-2}$ \cite{H}. When the radius is about $s$, the neighbors of $A_o$ are among the roughly $s^{d-4}$ trees meeting a box of radius $2s$. The contribution of each dyadic scale is at most a constant times $s^{d-4}\cdot s^{-2}=s^{d-6}$. This is summable exactly when $d=5$, which gives transience. In dimension six this bound is of order one per scale, which is the source of the logarithms in the six-dimensional return estimates. Recurrence for $d\ge6$ is proved separately, by tracing the walk onto the ray.

\emph{Repeated contacts and the threshold at dimension eight.} For $5\le d\le8$, any two trees of the forest are joined by a lattice edge \cite{BKPS}; more generally the trees form a component graph of diameter $\lceil(d-4)/4\rceil$. The scaling behind this also counts repeated contacts. Given one contact between two trees at the origin, a second contact at $w$ requires both trees to connect $0$ to $w$, with probability about $|w|^{2(4-d)}$. The expected number of contacts in $Q_r$ is therefore of order $\sum_{1\le|w|\le r}|w|^{8-2d}$, which grows like $\log r$ exactly when $d=8$ and stays bounded for $d\ge9$. For $d\ge8$ we bound the box resistance by a sum over contacts between $A_o$ and its neighbors, each weighted by the reciprocal of the number of contacts with the same neighbor. Since $\PP^{o,\mathrm w}(x\in A_o)\asymp(1+|x|)^{2-d}$ and $\sum_{x\in Q_R}(1+|x|)^{2-d}\asymp R^2$, bounded multiplicities would only give a bound of order $R^2$. In dimension eight, the logarithmic multiplicity lowers the weights by a factor of order $1/\log R$ (Proposition~\ref{prop:contact-multiplicity}). This produces the scales $R^2/\log R$ and $(n\log n)^{-2/3}$ for box resistance and return probability, respectively. The proof in Section~\ref{sec:contact-sharp} uses the local comparison methods of \cite{BKPS,HP}.

\subsection{Proof strategy}
\begin{figure}[ht]
\centering
\begin{tikzpicture}[x=0.772cm,y=0.78cm,font=\small,
  tr/.style={treenavy,line width=0.8pt,line cap=round,line join=round},
  ray/.style={rayteal,line width=1.2pt,line cap=round},
  ex/.style={extraorange,dashed,line width=0.8pt},
  v/.style={circle,fill=treenavy,inner sep=1.1pt}]
 \begin{scope}
  \node at (2.6,4.3) {Ordinary forest $F$};
  \fill[shadeA,rounded corners=6pt] (-0.7,-0.9) rectangle (3.0,2.6);
  \node[rayteal] at (0.2,2.25) {$K_e$};
  \draw[tr] (0.5,0.3)--(0.2,1.3)--(-0.2,2.0); \draw[tr] (0.2,1.3)--(0.9,1.9);
  \draw[tr] (0.5,0.3)--(1.5,0.7)--(2.2,1.7); \draw[tr] (1.5,0.7)--(1.9,-0.4);
  \draw[ray] (0.5,0.3)--(2.4,0.2);
  \draw[ray,extraorange,line width=1.6pt] (2.4,0.2)--(3.6,0.0);
  \node[extraorange!80!black,below] at (3.0,0.05) {$e$};
  \draw[ray,-{Stealth[length=5pt]}] (3.6,0.0)--(5.6,-0.3);
  \node[right] at (5.6,-0.3) {$\infty$};
  \draw[tr] (3.6,0.0)--(4.2,1.2)--(4.0,2.2); \draw[tr] (4.6,-0.15)--(4.9,-1.3);
  \foreach \p in {(0.2,1.3),(-0.2,2.0),(0.9,1.9),(1.5,0.7),(2.2,1.7),(1.9,-0.4),(2.4,0.2),(3.6,0.0),(4.2,1.2),(4.0,2.2),(4.6,-0.15),(4.9,-1.3)} \node[v] at \p {};
  \node[circle,fill=extraorange,inner sep=1.8pt] at (0.5,0.3) {};
  \node[below left] at (0.5,0.3) {$o$};
 \end{scope}
 \begin{scope}[xshift=8.2cm]
  \node at (2.6,4.3) {Wired forest $F^o=F-e$};
  \fill[extraorange!18,rounded corners=6pt] (-0.7,-0.9) rectangle (3.0,2.6);
  \node[extraorange!80!black] at (0.2,2.25) {$A_o$};
  \draw[tr] (0.5,0.3)--(0.2,1.3)--(-0.2,2.0); \draw[tr] (0.2,1.3)--(0.9,1.9);
  \draw[tr] (0.5,0.3)--(1.5,0.7)--(2.2,1.7); \draw[tr] (1.5,0.7)--(1.9,-0.4);
  \draw[tr] (0.5,0.3)--(2.4,0.2);
  \draw[tr,-{Stealth[length=4pt]}] (3.6,0.0)--(4.6,-0.15)--(6.0,-0.4); \draw[tr] (3.6,0.0)--(4.2,1.2)--(4.0,2.2); \draw[tr] (4.6,-0.15)--(4.9,-1.3);
  \draw[tr,-{Stealth[length=4pt]}] (2.9,3.2)--(4.2,3.4)--(6.0,3.6); \draw[tr] (1.0,3.2)--(2.9,3.2);
  \draw[tr,-{Stealth[length=4pt]}] (-1.3,-1.4)--(-1.5,-2.2)--(-2.2,-3.0); \draw[tr] (-1.3,-1.4)--(0.3,-1.6);
  \node[treenavy] at (5.6,0.5) {$K_1$};
  \node[treenavy] at (5.6,3.1) {$K_2$};
  \node[treenavy] at (-1.4,-3.0) {$K_3$};
  \draw[ex] (2.4,0.2)--(3.6,0.0); \draw[ex] (2.2,1.7)--(4.0,2.2); \draw[ex] (1.9,-0.4)--(4.9,-1.3);
  \draw[ex] (0.9,1.9)--(1.0,3.2); \draw[ex] (2.2,1.7)--(2.9,3.2);
  \draw[ex] (1.9,-0.4)--(0.3,-1.6);
  \foreach \p in {(0.2,1.3),(-0.2,2.0),(0.9,1.9),(1.5,0.7),(2.2,1.7),(1.9,-0.4),(2.4,0.2),(3.6,0.0),(4.2,1.2),(4.0,2.2),(4.6,-0.15),(4.9,-1.3),(1.0,3.2),(2.9,3.2),(4.2,3.4),(-1.3,-1.4),(-1.5,-2.2),(0.3,-1.6)} \node[v] at \p {};
  \node[circle,fill=extraorange,inner sep=1.8pt] at (0.5,0.3) {};
  \node[below left] at (0.5,0.3) {$o$};
  \node[align=left,anchor=west] at (6.4,1.6) {$m(A_o,K_1)=3$\\ $m(A_o,K_2)=2$\\ $m(A_o,K_3)=1$};
 \end{scope}
\end{tikzpicture}
\caption{The deletion identity of Lemma~\ref{prop:deletion}. Removing a ray edge $e$ leaves the finite side $K_e$. Summing over all such edges corresponds, up to the factor $r_G$, to summing over all boundary edges of $A_o$ under the wired law. Dashed edges are the lattice contacts between $A_o$ and its neighbors. Weighting each contact with $K$ by $1/m(A_o,K)$ counts $K$ once. Here $D_o=3$.}\label{fig:deletion}
\end{figure}
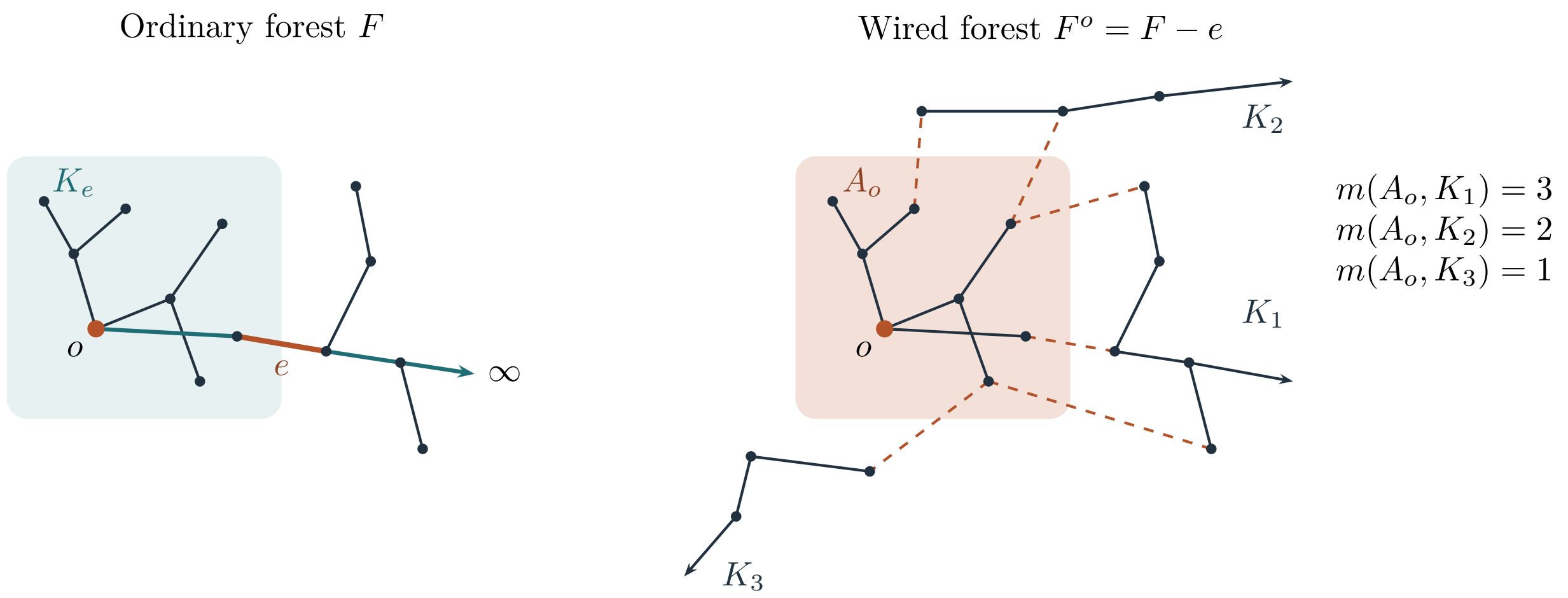

We use two representations of the forest. First, conditional on $H$, the tree $T$ has the wired uniform spanning-forest law $\WUSF(H)$ by the resampling theorem of Lyons, Peres and Sun \cite[Theorem~1.4]{LPS}. Averaging the signed current along its ray produces a current on $H$. A finite-cut identity expresses its voltage through reciprocal boundary sizes. Deleting the selected ray edge then changes the forest law to the law with the origin wired to infinity, denoted by $\PP^{o,\mathrm w}$ with expectation $\E^{o,\mathrm w}$. Under this law, reciprocal contact multiplicities count the components adjacent to the finite origin tree. If $D_o$ is their number and $\rho$ is the lattice radius of that tree about $o$, the conditional component bound and radius tail give
\[
 \E^{o,\mathrm w}[D_o\ind_{\{s\le\rho<2s\}}]
 \le Cs^{d-4}\PP^{o,\mathrm w}(\rho\ge s)\le Cs^{d-6}.
\]
These dyadic shell bounds are summable in dimension five; in dimension six each scale costs $O(1)$, leading to the logarithmic upper bound for box resistance (Figure~\ref{fig:deletion}). Section~\ref{sec:recurrence} makes the transience argument precise. This is a signed-current version of the random-path method in \cite[Section~2.5]{LP}.

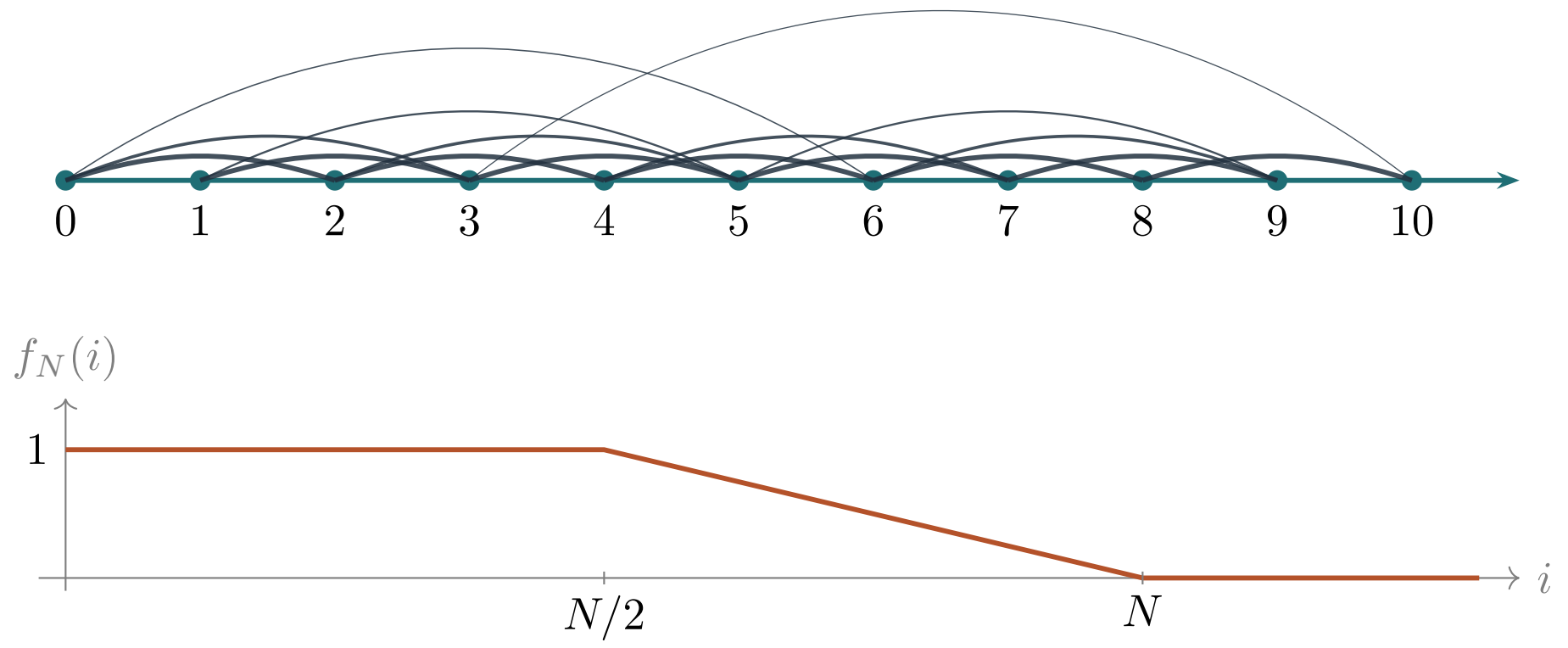
\begin{figure}[ht]
\centering
\begin{tikzpicture}[x=1.05cm,y=1cm,font=\small,
  v/.style={circle,fill=rayteal,inner sep=1.6pt}]
 \foreach \i in {0,...,10} {\node[v] (r\i) at (\i,0) {}; \node[below=2pt] at (\i,0) {$\i$};}
 \draw[rayteal,line width=1.1pt] (0,0)--(10,0);
 \draw[rayteal,line width=1.1pt,-{Stealth[length=5pt]}] (10,0)--(10.8,0);
 \foreach \i/\j/\w in {0/2/1.1,1/3/1.1,2/4/1.1,3/5/1.1,4/6/1.1,5/7/1.1,6/8/1.1,7/9/1.1,8/10/1.1,0/3/0.7,2/5/0.7,4/7/0.7,6/9/0.7,1/5/0.45,5/9/0.45,0/6/0.3,3/10/0.25}
  {\draw[treenavy,line width=\w pt,opacity=0.85] (\i,0) to[bend left={10+4*(\j-\i)}] (\j,0);}
 \begin{scope}[yshift=-3.1cm]
  \draw[latgray!60!black,->] (-0.2,0)--(10.8,0) node[right] {$i$};
  \draw[latgray!60!black,->] (0,-0.1)--(0,1.4) node[above] {$f_N(i)$};
  \draw[extraorange,line width=1.1pt] (0,1)--(4,1)--(8,0)--(10.5,0);
  \node[left] at (0,1) {$1$};
  \node[below] at (4,0) {$N/2$}; \node[below] at (8,0) {$N$};
  \draw[latgray!60!black] (4,0.05)--(4,-0.05); \draw[latgray!60!black] (8,0.05)--(8,-0.05);
 \end{scope}
\end{tikzpicture}
\caption{The walk on $H$ traced on the ray $\gamma$ (Lemma~\ref{lem:ray-trace}) is a reversible walk on $\{0,1,2,\ldots\}$ with conductances $a(i,j)$ between every pair of indices; their means are at most $C|i-j|^{1-d/2}$ by Lemma~\ref{sec:loop}. Extending a profile such as $f_N$ harmonically from the ray gives a voltage on $H$ of energy $\sum_{i<j}a(i,j)(f_N(i)-f_N(j))^2$, which tends to zero in mean for $d\ge7$, and for $d=6$ with a logarithmic profile.}\label{fig:trace}
\end{figure}

Second, tracing the walk onto $\gamma$ identifies conductances with conditional attachment probabilities through Wilson's algorithm \cite{Wilson,BLPS}. An injective loop count bounds their means by $C|i-j|^{1-d/2}$ (Figure~\ref{fig:trace}), and harmonic extension gives recurrence for $d\ge6$ and return lower bounds through a mass-transport confinement argument \cite{Lee23,LeeScaling}. For return upper bounds, applying the finite-cut identity to the charge $h-P_Hh$ of a smoothed function preserves signed cancellation, avoiding an extra logarithm in dimension six.

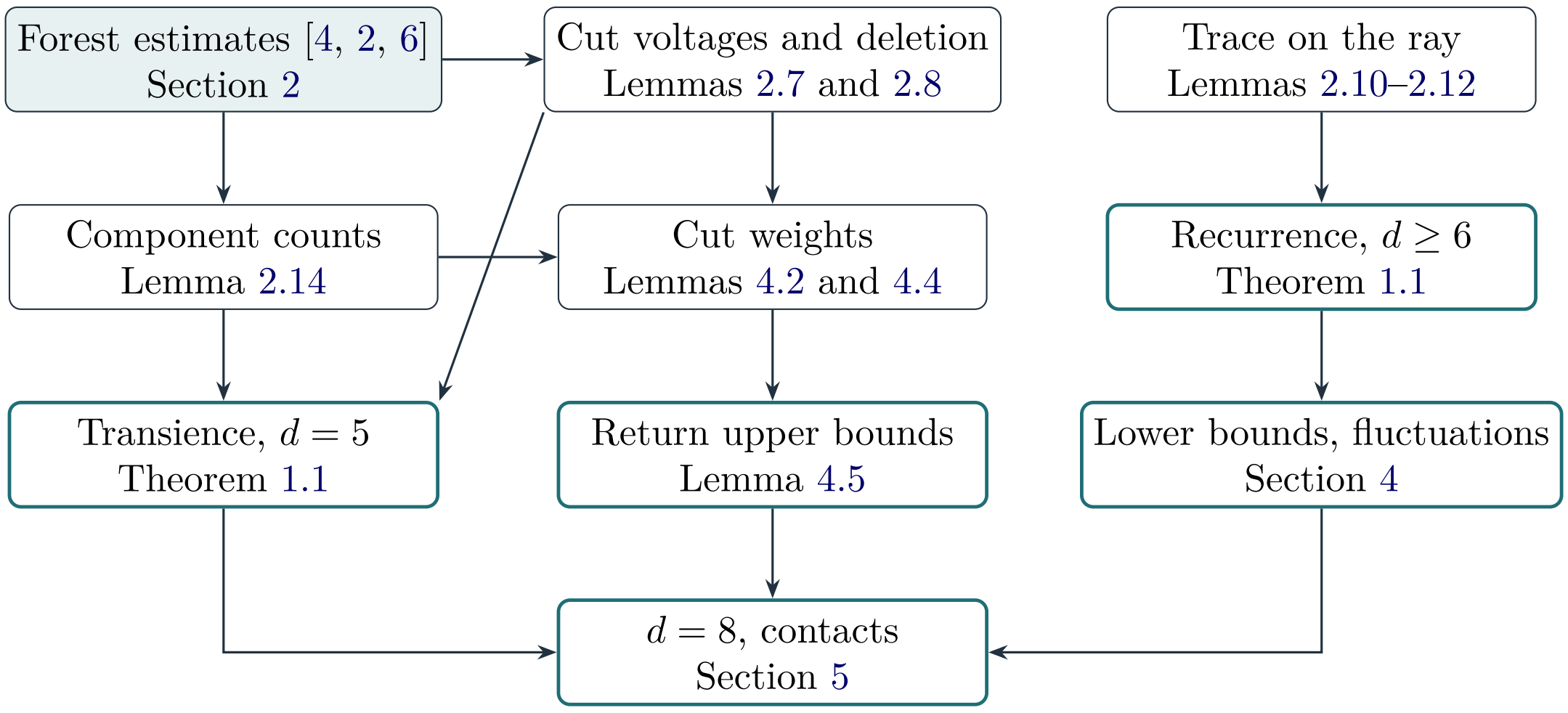
\begin{figure}[ht]
\centering
\begin{tikzpicture}[font=\small,x=1cm,y=1cm,
  box/.style={draw=treenavy,rounded corners=3pt,align=center,inner sep=3pt,fill=white,minimum width=3.9cm,minimum height=0.95cm},
  inp/.style={box,fill=shadeA},
  res/.style={box,draw=rayteal,line width=0.9pt},
  arr/.style={-{Stealth[length=5pt]},treenavy,line width=0.6pt}]
 \node[inp] (bkps) at (0,0) {Forest estimates \cite{BKPS,BJ,H}\\ Section~\ref{sec:prelim}};
 \node[box] (del) at (5,0) {Cut voltages and deletion\\ Lemmas~\ref{prop:finite-cut} and~\ref{prop:deletion}};
 \node[box] (trace) at (10,0) {Trace on the ray\\ Lemmas~\ref{lem:ray-trace}--\ref{sec:cutoffs}};
 \node[box] (pack) at (0,-1.8) {Component counts\\ Lemma~\ref{lem:packing}};
 \node[box] (cuts) at (5,-1.8) {Cut weights\\ Lemmas~\ref{lem:dipoles} and~\ref{lem:cut-statistics}};
 \node[res] (rec) at (10,-1.8) {Recurrence, $d\ge6$\\ Theorem~\ref{thm:main}};
 \node[res] (trans) at (0,-3.6) {Transience, $d=5$\\ Theorem~\ref{thm:main}};
 \node[res] (upper) at (5,-3.6) {Return upper bounds\\ Lemma~\ref{lem:returnupper}};
 \node[res] (lower) at (10,-3.6) {Lower bounds, fluctuations\\ Section~\ref{sec:heat}};
 \node[res] (d8) at (5,-5.4) {$d=8$, contacts\\ Section~\ref{sec:contact-sharp}};
 \draw[arr] (bkps)--(pack);
 \draw[arr] (bkps)--(del);
 \draw[arr] (pack)--(trans);
 \draw[arr] (pack)--(cuts);
 \draw[arr] (del)--(cuts);
 \draw[arr] (del.south west)--(trans.north east);
 \draw[arr] (cuts)--(upper);
 \draw[arr] (trace)--(rec);
 \draw[arr] (rec)--(lower);
 \draw[arr] (upper)--(d8);
 \draw[arr] (lower)|-(d8);
 \draw[arr] (trans)|-(d8);
\end{tikzpicture}
\caption{Structure of the proofs. Shaded: inputs from the literature. Outlined in teal: the results of Theorems~\ref{thm:main} and~\ref{thm:returns}. The eight-dimensional bounds also use the comparison methods of \cite{BKPS,HP} and the cut formula of Lemma~\ref{lem:terminal}.}\label{fig:structure}
\end{figure}

The fluctuations in dimensions five through seven come from tilting a middle segment of the first Wilson walk, which makes the ray locally ballistic and lowers the voltage energy at a bounded likelihood cost. In dimension eight, repeated contacts between neighboring components reduce reciprocal cut weights by a logarithmic factor. The contact estimates build on \cite{BKPS,HP} and yield the sharp box-resistance and mean-return scales. The passage from volume and resistance estimates to return probabilities is related to the general methods of \cite{BJKS,KM}; here the cut and trace representations provide estimates adapted to the induced graph. Figure~\ref{fig:structure} shows how the parts depend on each other.

\subsection{Organization of the paper}
Section~\ref{sec:prelim} supplies the forest estimates and representations, ending with the component counts of Lemma~\ref{lem:packing}. Section~\ref{sec:recurrence} proves Theorem~\ref{thm:main}, and Section~\ref{sec:terminal} gives a cut formula for box resistance and its sharp six- and seven-dimensional scales. Section~\ref{sec:heat} proves Theorem~\ref{thm:returns} apart from the eight-dimensional row: lower bounds in Section~\ref{sec:heatlower}, upper bounds in Section~\ref{sec:heatupper}, fluctuations in Section~\ref{sec:heat6}, and the remaining six-dimensional estimates in Section~\ref{sec:d6}. Section~\ref{sec:contact-sharp} proves the eight-dimensional bounds, and Section~\ref{sec:questions} lists open questions.

\medskip

\noindent {\bf Acknowledgements:} We thank Xinxin Chen for useful discussions on this work. XL and RL are supported by National Natural Science Foundation of China (No.~12631006), the National Key R\&D Program of China (No.~2021YFA1002700) and Beijing Natural Science Foundation (JQ26001). Yuval Peres acknowledges support from the National Natural Science Foundation of China under grant~\mbox{RFIS-W2531011}.

\medskip

\noindent {\bf AI Statement:} The development of the key idea in this work benefited from discussions with commercial large language models. The authors take full responsibility for the correctness of this work.

\section{Preliminaries}\label{sec:prelim}
We collect forest estimates, derive cut and deletion identities, and construct localized voltages through the trace on the ray. The final lemma counts the components meeting a box, also under the conditioning produced by deletion.

\subsection{Notation}\label{sec:notation}
We collect general conventions here. The principal notation is indexed in Table~\ref{tab:notation}.

\paragraph{General conventions.}
Except in Remark~\ref{rem:conductances}, all edge conductances are one. For $x\in\mathbb R^d$, write $|x|_1$ and $|x|_2$ for the $\ell^1$ and Euclidean norms, and set
\[
 |x|=\|x\|_\infty,\qquad \langle x\rangle=1+\|x\|_\infty.
\]
Write $e_1,\ldots,e_d$ for the coordinate unit vectors and $I_d$ for the $d\times d$ identity matrix. For a vertex set $A$, $|A|$ is its cardinality and $\diam_{\rm ext}A=\sup_{x,y\in A}|x-y|$. Constants $c,C$ are positive and deterministic, may depend on $d$, and may change between occurrences; further fixed-parameter dependence is indicated by subscripts or stated locally. We use $a\vee b=\max\{a,b\}$, $a\wedge b=\min\{a,b\}$, $a_+=a\vee0$, and $\mathbb N_0=\{0,1,\ldots\}$. The symbols $\ind_E$ and $\mathcal L(Y)$ denote an event or set indicator and the law of $Y$; $\ind_x$ is the indicator of the singleton $\{x\}$, viewed as a function on vertices. For probability laws, total variation distance is $d_{\rm TV}(\mu,\nu)=\sup_E|\mu(E)-\nu(E)|$.

\paragraph{Graphs and walks.}
For a graph $J$, let $J[A]$ be the subgraph induced by $A$, write $x\sim_Jy$ for adjacency, and let $E_J(A,B)$ be the set of edges with one endpoint in $A$ and the other in $B$. Thus $\partial_J A=E_J(A,V(J)\setminus A)$ is the edge boundary and $b_J(A)=|\partial_J A|$. For disjoint lattice vertex sets, $A\sim B$ means $E_G(A,B)\ne\varnothing$. For a vertex function $\mu$ and a finite set $A$, write $\mu(A)=\sum_{x\in A}\mu(x)$, in particular $\pi(A)=\sum_{x\in A}\pi(x)$.

For simple random walk on $J$, use $P_x^J,E_x^J$ for its law and expectation, $p_n^J(x,y)=P_x^J(X_n=y)$ for its kernel, $P_Jf(x)=E_x^J[f(X_1)]$ for its transition operator, and $\pi_J(x)=\deg_J(x)$ for its reversible measure. We omit $J$ from $P_x^J,E_x^J$ only when the underlying graph is specified; lattice kernels retain the superscript $G$, whereas $p_n$ denotes $p_n^H(o,o)$. For a walk on a specified graph, put
\[
 \tau_A=\inf\{n\ge0:X_n\in A\},\qquad
 \tau_A^+=\inf\{n\ge1:X_n\in A\},
\]
with $\inf\varnothing=\infty$ and $\tau_x=\tau_{\{x\}}$, $\tau_x^+=\tau_{\{x\}}^+$. For real $f,g\in\ell^2(V(J),\pi_J)$, the inner product and norm are
\[
 \langle f,g\rangle_{\pi_J}=\sum_{x\in J}\pi_J(x)f(x)g(x),\qquad
 \|f\|_{\pi_J}^2=\langle f,f\rangle_{\pi_J}.
\]
For any real vertex function $f$, its energy is
\[
 \mathcal E_J(f)=\sum_{\{x,y\}\in E(J)}(f(x)-f(y))^2.
\]
For $f\in\ell^2(V(J),\pi_J)$ this equals $\langle f,(I-P_J)f\rangle_{\pi_J}$, with $I$ the identity. For $J=H$, we abbreviate $\pi_H$ to $\pi$.
On the infinite graph $H$, resistances are understood as follows. For $B\subset V(H)$ with finite complement, $\res^H(o,B)$ is the resistance of the finite network obtained by identifying $B$ to a single vertex. The escape resistance is the wired one: $\res^H(o,\infty)^{-1}=\inf\{\mathcal E_H(f):f(o)=1,\ f\text{ finitely supported}\}$, while $\res^H(o,\infty)$ itself is the infimum of the energies of unit flows from $o$ to infinity by Thomson's principle \cite[Chapters~2 and~9]{LP}. For two vertices, $\res^H(o,x)$ is the wired resistance, the infimum of the energies of unit flows from $o$ to $x$. On $H$ this also equals the free resistance: by \cite[Theorem~1.4 and Corollary~1.6]{LPS}, the free and wired forest laws on $H$ coincide almost surely, so every harmonic function of finite Dirichlet energy on $H$ is constant \cite[Chapters~9--10]{LP}.

\paragraph{Paths and forest laws.}
For a path $\eta$, $|\eta|$ is its number of edges, $\eta[a,b]$ is its indexed segment, and $\eta\oplus\xi$ denotes concatenation with the common endpoint counted once. A path is regarded as its vertex set in intersections and as its edge set in edge sums. Write $\LE(\eta)$ for chronological loop erasure, $[x,y]_T$ for the unique path in a tree $T$, and $\Ray_x(F)$ for the infinite ray from $x$ in a one-ended forest. Each such tree is oriented toward its unique end, so the outgoing edge at a vertex is its first ray edge. We write $x\leftrightarrow y$ when $x,y$ belong to the same component of $F$. For a connected graph $J$, $\UST(J)$ denotes the uniform spanning-tree law when $J$ is finite, and $\WUSF(J)$ the wired uniform spanning-forest law when $J$ is infinite.

The $x$-wired law $\PP^{x,\mathrm w}$, with expectation $\E^{x,\mathrm w}$, is obtained by identifying $x$ with the sink in finite wired boxes, and then separating them while retaining the original lattice edges. Let $F^x$ have this law and write $A_x$ for its finite component containing $x$. All other components are infinite.

\paragraph{Principal notation.}
Table~\ref{tab:notation} lists selected notation for the main objects, laws and estimates.

\begin{longtable}{@{}>{\raggedright\arraybackslash}p{0.24\textwidth}>{\raggedright\arraybackslash}p{0.53\textwidth}>{\raggedright\arraybackslash}p{0.16\textwidth}@{}}
\caption{Principal notation.}\label{tab:notation}\\
\toprule
\textbf{Symbol} & \textbf{Meaning} & \textbf{Defined in}\\
\midrule
\endfirsthead
\toprule
\textbf{Symbol} & \textbf{Meaning} & \textbf{Defined in}\\
\midrule
\endhead
\midrule
\endfoot
\bottomrule
\endlastfoot
$G$, $F$, $T_x$, $H_x$ & Lattice, its WUSF, the tree containing $x$, and its induced graph; $T=T_o$, $H=H_o$. & \S\ref{sec:model}\\
$\PP$, $\E$; $P_x^H$, $E_x^H$ & Forest law and expectation; walk law and expectation on the fixed graph $H$. & \S\ref{sec:model}\\
$\PP^{x,\mathrm w}$, $\E^{x,\mathrm w}$, $F^x$, $A_x$ & $x$-wired law, expectation, forest and finite component containing $x$. & \S\ref{sec:notation}\\
$P_H$, $\mathbf Q_H$ & Transition operator of the walk on $H$, and the conditional tree law given $H$. & \eqref{eq:resampling}\\
$p_n$ & Return probability $p_n^H(o,o)$. & \S\ref{sec:results}\\
$Q_R(x)$, $d_J$, $B_J(x,R)$ & Lattice cube, graph distance and closed graph ball; $Q_R=Q_R(o)$. & \S\ref{sec:results}\\
$\mathcal R_R$, $\mathcal R_\infty$ & Resistance to the exterior of $Q_R$ and to infinity. & \S\ref{sec:results}\\
$\gamma=(r_i)_{i\ge0}$, $\mathcal B_i$ & Ray from $o$ and finite bush attached at $r_i$. & \S\ref{sec:model}, \S\ref{sec:cuts}\\
$\pi$, $\mathcal E_H$ & Reversible degree measure and Dirichlet energy on $H$. & \S\ref{sec:notation}\\
$g_G(x,y)$ & Lattice Green function. & Lemma~\ref{lem:lattice-walk}\\
$r_G$, $D_o$, $m(A,K)$ & Lattice escape resistance, number of components adjacent to $A_o$, and number of lattice edges between $A,K$. & Lemma~\ref{prop:deletion}\\
$a(i,j)$ & Conductance of the walk on $H$ traced onto the ray. & \eqref{eq:trace}\\
$k_d(N)$ & Common energy scale for localized cutoff voltages. & \eqref{energies}\\
$L_A$, $w_A$ & Event that $T[A]$ is connected with a sole forest exit, and conditional cut weight $\mathbf Q_H(L_A)/b_H(A)$. & Lemmas~\ref{lem:terminal}, \ref{lem:dipoles}\\
$\rho$, $\operatorname{Cap}_a(0)$ & Radius of $A_o$ about $o$ and capacity of the root for the walk traced on the ray. & Lemma~\ref{lem:packing}, \S\ref{sec:recurrence}\\
$p(x)$, $\mathcal D_x$, $N_r(x)$ & Ordered component pair across $(x,x+e)$, its distinctness event, and the number of matching pairs in $Q_r(x)$. & \S\ref{sec:contact-sharp}\\
$q(r)$, $\phi_{r,S}(x)$ & Mean reciprocal contact count with distinct components, and its version computed from a box. & \S\ref{sec:contact-sharp}, \S\ref{sec:contact-localization}\\
\end{longtable}

\subsection{Uniform spanning trees and forests}\label{sec:ust-background}
We record the following standard results on uniform spanning forests and the lattice random walk.

\begin{lemma}[{\cite[Theorem~5.1]{BLPS}; \cite[Exercise~10.8]{LP}; \cite[Theorem~1.4]{LPS}}]\label{lem:wilson-facts}
For $d\ge5$, conditional on exposed Wilson paths, the remaining edges have the forest law rooted at their union and infinity. Adding roots stochastically decreases the edge set, in finite wired graphs and their limits. Moreover,
\begin{equation}\label{eq:resampling}
 \mathcal L(T\mid H)=\WUSF(H)=:\mathbf Q_H.
\end{equation}
\end{lemma}
To see the conditional assertion in \eqref{eq:resampling}, resample a WUSF on each induced component. By \cite[Theorem~1.4]{LPS} the new forest has the original law. Its partition into components refines the original partition. For each pair of vertices the two connection events thus agree almost surely, since they have the same probability. Hence the partitions agree almost surely. The induced graph at $o$ can therefore be recovered from the resampled forest, and the joint law of this graph and its tree is the original joint law of $(H,T)$.

\begin{lemma}[{\cite[Propositions~2.1.2, 2.4.6 and Theorem~4.3.1]{LL}; \cite[Lemma~2.10]{BKPS}}]\label{lem:lattice-walk}
Let $p_n^G(x,y)$ be the simple random-walk kernel on $G=\Z^d$ and let $g_G(x,y)=\sum_{n\ge0}p_n^G(x,y)$. For $d\ge5$,
\begin{gather}
 p_n^G(x,y)\le C(n+1)^{-d/2}
   \exp\{-c|x-y|_2^2/(n+1)\},\label{eq:lattice-heat}\\
 g_G(x,y)\asymp(1+|x-y|)^{2-d},\quad
 \sum_z\!g_G(x,z)g_G(z,y)\le C(1+|x-y|)^{4-d}.\label{eq:lattice-green}
\end{gather}
For $x\ne y$ and a constant $\mathfrak a_d>0$,
\begin{equation}\label{eq:green-expansion}
 g_G(x,y)=\mathfrak a_d|x-y|_2^{2-d}+O(|x-y|_2^{-d}).
\end{equation}
If the lattice walk $S$ starts at $x$, then
$P_x^G(\max_{j\le n}|S_j-x|\ge R)\le C\exp(-cR^2/n)$ for $n,R\ge1$.
\end{lemma}

\begin{lemma}[{\cite[Theorems~4.2--4.3]{BKPS}; \cite[Theorem~5.4]{BJ}}]\label{lem:forest-volumes}
For $d\ge5$, $R\ge1$, $V_R=|H\cap Q_R|$, integers $k\ge1$, and $0<a\le1$,
\begin{gather}
 \PP(x\leftrightarrow y\text{ in }F)\asymp(1+|x-y|)^{4-d},\label{eq:two-point}\\
 \E[V_R^k]\le C_{k,d}R^{4k},\quad
 \PP(V_R\le aR^4)\le Ce^{-ca^{-1/3}}.\label{volumes}
\end{gather}
Consequently $V_R=R^{4+o(1)}$ almost surely.
\end{lemma}
For the lower tail in \eqref{volumes}, wire the complement of $Q_{4R}$ to a vertex and remove that vertex from the sampled tree. By root monotonicity its internal edge set can be coupled below $F$ (first use a larger wired box and then take a limit). Its origin component inside $Q_R$ is therefore contained in $H\cap Q_R$, so the finite-volume estimate of \cite[Theorem~5.4]{BJ} applies directly.

\begin{lemma}[{\cite[Theorem~1.1, Proposition~3.12 and Lemma~5.8]{BJ}}]\label{lem:intrinsic-volume}\label{lem:ray-length}
For $d\ge5$ and $m,R,\lambda\ge1$,
\begin{equation}\label{eq:intrinsic-volume-tails}
 \PP(|B_T(o,m)|\ge\lambda m^2)\le Ce^{-c\lambda},
 \quad
 \PP(|B_T(o,m)|\le\lambda^{-1}m^2)\le Ce^{-c\lambda^{1/5}}.
\end{equation}
Thus $|B_T(o,m)|=m^{2+o(1)}$ almost surely, and every fixed positive and negative moment of $|B_T(o,m)|/m^2$ is bounded uniformly in $m$. If
$\sigma_R=\min\{j:r_j\notin Q_R\}$, then
\begin{equation}\label{eq:ray-exit-tails}
 \PP(\sigma_R\le\lambda^{-1}R^2)\le Ce^{-c\lambda},
 \quad
 \PP(\sigma_R\ge\lambda R^2)\le Ce^{-c\lambda}.
\end{equation}
\end{lemma}

\begin{lemma}[Past and wired-root estimates]\label{lem:past-estimates}
For $d\ge5$, let $\Past_F(x)$ be the finite side of the outgoing edge at $x$, including $x$, and write $\operatorname{rad}_{A_x}(x)=\max_{y\in A_x}d_{A_x}(x,y)$. For $R,m,t\ge1$,
\begin{gather}
 \PP(\diam_{\rm ext}\Past_F(x)\ge R)\le CR^{-2},\quad
 \PP^{x,\mathrm w}(\diam_{\rm ext}A_x\ge R)\le CR^{-2},\label{eq:Htails}\\
 \PP^{x,\mathrm w}(\operatorname{rad}_{A_x}(x)\ge m)\le C/m,\quad
 \PP^{x,\mathrm w}(|A_x|\ge t)\le C/\sqrt t,\label{eq:wired-past-tails}\\
 \E^{x,\mathrm w}[|B_{A_x}(x,m)|]\le C(m+1).\label{eq:wired-truncated-volume}
\end{gather}
The upper bounds for $A_x$ also hold with $A_x$ replaced by $\Past_F(x)$.
\end{lemma}
The two bounds in \eqref{eq:Htails} follow from \cite[Theorems~1.4 and~7.2]{H}, respectively. Equation~\eqref{eq:wired-past-tails} follows from \cite[Theorem~7.1]{H}, and \eqref{eq:wired-truncated-volume} is \cite[Proposition~6.2]{H} with moment order one.

\begin{lemma}[Conditional bush domination {\cite[Lemma~2.1]{H}}]\label{lem:bush-domination}
Given the futures of finitely many Wilson starting vertices, delete the exposed future edges. For each vertex $v$ of these futures, the remaining branch attached at $v$ is stochastically dominated, by edge inclusion, by $A_v$ under $\PP^{v,\mathrm w}$. In particular this holds for the bush at $r_j$ conditional on the entire ray $\gamma$. The same bound holds in the presence of additional Wilson roots.
\end{lemma}
Indeed, the unexposed edges have the law rooted on all the exposed vertices and at infinity. Keeping only $v$ and infinity as roots gives the domination by Lemma~\ref{lem:wilson-facts}. The two trees in this inclusion coupling have the same distances from $v$ on their common branch, so the radius and intrinsic-ball estimates apply there as well.

\subsection{Finite cuts and deletion}\label{sec:cuts}
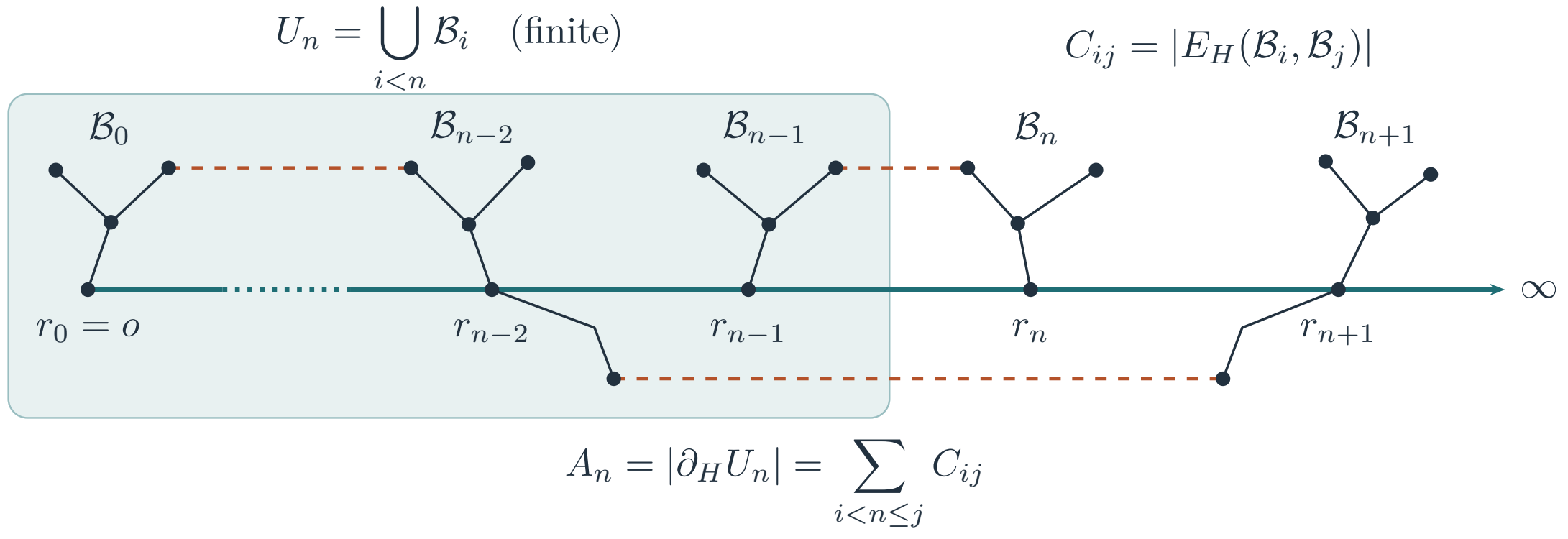
\begin{figure}[htbp]
\centering
\begin{tikzpicture}[x=1.18cm,y=1cm,font=\small,
  tree/.style={draw=treenavy,line width=.65pt,line cap=round,line join=round},
  ray/.style={draw=rayteal,line width=1.25pt},
  extra/.style={draw=extraorange,dashed,line width=.85pt},
  vertex/.style={circle,fill=treenavy,inner sep=1.35pt},
  every node/.style={text=treenavy}]
  \fill[shadeA,rounded corners=5pt] (-.62,-1.18) rectangle (6.25,1.80);
  \draw[rayteal!45,rounded corners=5pt,line width=.45pt]
    (-.62,-1.18) rectangle (6.25,1.80);
  \node at (2.82,2.22) {$U_n=\displaystyle\bigcup_{i<n}\mathcal B_i\quad\text{(finite)}$};
  \node at (8.82,2.22) {$C_{ij}=|E_H(\mathcal B_i,\mathcal B_j)|$};
  \draw[ray] (0,0)--(1.05,0);
  \draw[ray,dotted] (1.05,0)--(2.05,0);
  \draw[ray] (2.05,0)--(3.15,0)--(5.15,0)--(7.35,0)--(9.75,0);
  \draw[ray,-{Stealth[length=4pt]}] (9.75,0)--(11.05,0);
  \node[right] at (11.05,0) {$\infty$};
  \draw[tree] (0,0)--(.18,.62)--(-.25,1.10);
  \draw[tree] (.18,.62)--(.63,1.12);
  \draw[tree] (3.15,0)--(2.97,.60)--(2.52,1.12);
  \draw[tree] (2.97,.60)--(3.43,1.17);
  \draw[tree] (3.15,0)--(3.95,-.35)--(4.10,-.82);
  \draw[tree] (5.15,0)--(5.31,.60)--(5.83,1.12);
  \draw[tree] (5.31,.60)--(4.80,1.10);
  \draw[tree] (7.35,0)--(7.25,.61)--(6.86,1.12);
  \draw[tree] (7.25,.61)--(7.86,1.10);
  \draw[tree] (9.75,0)--(10.02,.66)--(9.65,1.18);
  \draw[tree] (10.02,.66)--(10.47,1.06);
  \draw[tree] (9.75,0)--(9.00,-.35)--(8.85,-.82);
  \draw[extra] (.63,1.12)--(2.52,1.12);
  \draw[extra] (5.83,1.12)--(6.86,1.12);
  \draw[extra] (4.10,-.82)--(8.85,-.82);
  \foreach \p in {(0,0),(3.15,0),(5.15,0),(7.35,0),(9.75,0),
    (.18,.62),(-.25,1.10),(.63,1.12),(2.97,.60),(2.52,1.12),
    (3.43,1.17),(4.10,-.82),(5.31,.60),(5.83,1.12),(4.80,1.10),
    (7.25,.61),(6.86,1.12),(7.86,1.10),(10.02,.66),(9.65,1.18),
    (10.47,1.06),(8.85,-.82)} \node[vertex] at \p {};
  \node[below=4pt] at (0,0) {$r_0=o$};
  \node[below=4pt] at (3.15,0) {$r_{n-2}$};
  \node[below=4pt] at (5.15,0) {$r_{n-1}$};
  \node[below=4pt] at (7.35,0) {$r_n$};
  \node[below=4pt] at (9.75,0) {$r_{n+1}$};
  \node at (.17,1.48) {$\mathcal B_0$};
  \node at (3.00,1.48) {$\mathcal B_{n-2}$};
  \node at (5.28,1.48) {$\mathcal B_{n-1}$};
  \node at (7.40,1.48) {$\mathcal B_n$};
  \node at (10.04,1.48) {$\mathcal B_{n+1}$};
  \node at (5.36,-1.78)
    {$A_n=|\partial_HU_n|=\displaystyle\sum_{i<n\le j}C_{ij}$};
\end{tikzpicture}
\caption{The shaded finite set $U_n$ and its boundary. Solid edges belong to $T$, dashed edges to $H\setminus T$. The cut size $A_n$ includes all crossing edges, including $(r_{n-1},r_n)$.}
\label{fig:bush-cuts}
\end{figure}

For $i\ge0$, let $\mathcal B_i$ be the finite branch at $r_i$ obtained by removing the edges of $\gamma$; we also call it the bush at $r_i$. For $n\ge1$ and $i,j\ge0$, define
\begin{equation}\label{eq:Vn}
 U_n=\bigcup_{i<n}\mathcal B_i=\Past_F(r_{n-1}),\qquad
 C_{ij}=|E_H(\mathcal B_i,\mathcal B_j)|,\qquad
 A_n=\sum_{i<n\le j}C_{ij}.
\end{equation}
Thus $A_n=b_H(U_n)$; Figure~\ref{fig:bush-cuts} illustrates the cut. The next lemma expresses resistance through these cut sizes.

\begin{lemma}[Cut voltages]\label{prop:finite-cut}\label{thm:infinite-cut}
For the general-graph statements in this lemma and their proofs, $G$ is a connected, locally finite graph and $\PP,\E$ refer to a tree $T$ sampled from $\UST(G)$ when $G$ is finite and from $\WUSF(G)$ otherwise. For finite $G$, vertices $a,b\in V(G)$, and $e\in[a,b]_T$, let $K_e$ be the component of $T-e$ containing $a$. Then
\begin{equation}\label{eq:finite-cut}
 \res^G(a,b)=\E_{\UST(G)}\sum_{e\in[a,b]_T}|\partial_GK_e|^{-1}.
\end{equation}
If $G$ is infinite and its WUSF is a single one-ended spanning tree, fix $o\in V(G)$ and, for $e\in\Ray_o(T)$, let $K_e$ be the finite component of $T-e$ containing $o$. Then
\begin{equation}\label{eq:SG}
 S_G=\sum_{e\in\Ray_o(T)}|\partial_GK_e|^{-1},\quad
 \E S_G=\res^G(o,\infty),
\end{equation}
with equality in $[0,\infty]$. When this value is finite, the unit-source voltage is
\begin{equation}\label{eq:u}
 u(x)=\E\sum_{e\in\Ray_o(T)}\frac{\ind_{\{x\in K_e\}}}{|\partial_GK_e|}.
\end{equation}
It vanishes at infinity, and its gradient is the mean signed current along the ray.
For the original lattice forest and its induced component, conditioning on $H$ gives
\begin{equation}\label{eq:conditional-resistance}
 \res^H(o,\infty)=\E\left[\left.\sum_{n\ge1}A_n^{-1}\right|H\right].
\end{equation}
\end{lemma}
\begin{proof}
The two-component forest formula \cite[Theorem~3]{KKW} gives $\res^G(a,b)=Z_2(a,b)/Z_1$, where $Z_1$ counts spanning trees and $Z_2(a,b)$ separating two-component forests. Each separating forest has $|\partial_GK|$ tree completions, so the reciprocal weights in \eqref{eq:finite-cut} count it once.

For the infinite identity, couple the trees of a wired exhaustion so that their edges converge locally to $T$ almost surely. For each finite side $K_e$, the internal tree and every edge of $\partial_GK_e$ eventually agree with the limit, so the same side occurs on the finite path to the sink. Fatou's lemma gives
\begin{equation}\label{eq:fatou}
 \E S_G\le\res^G(o,\infty).
\end{equation}
This suffices if the left side is infinite. Otherwise, for finite connected $K\ni o$, let $L_K^e$ be the event that $T[K]$ is connected with sole exit $e\in\partial_GK$, and put $L_K=\bigcup_eL_K^e$. Finite-volume edge exchange passes to these cylinder events and gives
\begin{equation}\label{eq:cylinder-exchange}
 \PP(L_K^e)=\PP(L_K^f)=\PP(L_K)/|\partial_GK|.
\end{equation}
Summing over $K$ containing $x$ but not its neighbor $y$ identifies the corresponding part of \eqref{eq:u} with the probability that the ray traverses $(x,y)$. Hence
\begin{equation}\label{eq:gradient}
 u(x)-u(y)=\PP((x,y)\in\Ray_o(T))-\PP((y,x)\in\Ray_o(T)),\qquad
 \sum_{y\sim_Gx}(u(x)-u(y))=\ind_{\{x=o\}}.
\end{equation}
The series at $o$ is summable, so $u$ tends to zero outside finite unions of its finite sets $K$. Iterating its Poisson equation gives
\begin{equation}\label{eq:green-bound}
 \deg_G(o)^{-1}\sum_{k<m}P_o^G(X_k=o)=u(o)-E_o^G[u(X_m)]\le u(o).
\end{equation}
Thus $\E S_G<\infty$ forces transience, and $\res^G(o,\infty)\le u(o)=\E S_G$. Together with \eqref{eq:fatou} this proves equality, including when $\res^G(o,\infty)=\infty$. Finally test the Poisson equation against $(u-t)_+$, which has finite support. Summation by parts and monotone convergence as $t\downarrow0$ give $\mathcal E_G(u)=u(o)$, identifying the finite-energy unit current. The conditional identity follows from Lemma~\ref{lem:wilson-facts} and $A_n=b_H(U_n)$.
\end{proof}
Returning to $G=\Z^d$, deleting the chosen ray edge changes the forest law; reciprocal contact counts cancel its possible choices.
\begin{lemma}[Deletion and neighboring components]\label{sec:deletion}\label{prop:deletion}\label{thm:environment}\label{cor:one-edge}
Write $r_G=\res^{\Z^d}(o,\infty)$ and let $D_o$ count the components of $F^o$ adjacent to $A_o$. For every nonnegative measurable function $\Phi$ of a forest and an edge,
\begin{gather}
 \E\sum_{e\in\Ray_o(F)}\Phi(F-e,e)
 =r_G\E^{o,\mathrm w}\sum_{e\in\partial_GA_o}\Phi(F^o,e),\label{eq:deletion}\\
 \E\res^H(o,\infty)=\E\sum_{n\ge1}A_n^{-1}=r_G\E^{o,\mathrm w}D_o.\label{eq:environment}
\end{gather}
Moreover, $F^o,F$ admit a coupling with
\begin{equation}\label{eq:one-edge}
 F^o\subseteq F,\quad |E(F)\setminus E(F^o)|=1\quad\hbox{almost surely}.
\end{equation}
\end{lemma}
\begin{proof}
In a finite wired box, deleting an $o$--sink path edge and adding an edge between the resulting components are inverse operations. The number of separating two-component forests divided by the tree count is the finite-volume resistance. Fix a finite connected $A\ni o$ and $e\in\partial_G A$. Let $L_A^e$ be the event that $F[A]$ is connected and $e$ is its sole forest exit. Both $L_A^e$ and $\{A_o=A\}$ depend only on edges in or incident to $A$. Testing the finite bijection against a bounded cylinder function $\Psi$ and taking the wired limit gives
\begin{equation}\label{eq:marked}
 \E[\ind_{L_A^e}\Psi(F-e)]=r_G\E^{o,\mathrm w}[\ind_{\{A_o=A\}}\Psi(F^o)].
\end{equation}
The right side is independent of $e$, as in the uniform-exit identity \eqref{eq:cylinder-exchange}. For each fixed $(A,e)$, the two finite measures in \eqref{eq:marked} agree on cylinder events, hence on the full forest $\sigma$-field. We may thus use any nonnegative measurable test. Summing over $(A,e)$ gives \eqref{eq:deletion} by Tonelli: every ray edge has a finite origin side, and $A_o$ is finite almost surely.

Let $K$ be the component containing the outer endpoint of $e$ in the deleted forest. Restoring $e$ gives the induced graph $G[A\cup K]$, with $m(A,K)=|E_G(A,K)|$ edges across $A$. This count is measurable since connectivity is the union of the events that a given finite path is present. Testing with $1/m(A,K)$ therefore gives $A_n^{-1}$ on the ordinary side and a total weight of one for each neighboring component on the wired side. Equation~\eqref{eq:conditional-resistance} now gives \eqref{eq:environment}.

In each finite wired box, root monotonicity gives an inclusion coupling (equivalently, use determinantal domination \cite[Theorem~6.2]{Lyons2003}). The edge counts differ by one. Inclusion and the condition of having at most one differing edge are closed in the product topology, since a violation is witnessed by at most two edges. They persist in any subsequential limit. As $A_o$ is finite and the ordinary origin tree is infinite, the extra edge must join $A_o$ to an infinite component. It is the unique exit from $A_o$ in $F$ and hence lies on $\Ray_o(F)$.
\end{proof}

\begin{remark}[Positive edge conductances]\label{rem:conductances}
Let $(G,c)$ be a connected, locally finite network with $0<c(e)<\infty$ and finite-volume spanning-tree weights $\prod_{e\in T}c(e)$. All forest laws in this remark refer to this network. Under the respective hypotheses of Lemma~\ref{prop:finite-cut}, replace $|\partial_GA|$ by $c(\partial_GA)=\sum_{e\in\partial_GA}c(e)$ in \eqref{eq:finite-cut}, \eqref{eq:SG} and \eqref{eq:u}, using the weighted forest laws and resistances. The voltage gradient carries current $c(x,y)(u(x)-u(y))$, where $c(x,y)$ is the conductance of the edge $\{x,y\}$.

For Lemma~\ref{prop:deletion}, assume in addition that $G$ is transient, its WUSF components are one-ended, and its $o$-wired component is finite almost surely. With $r_G=\res^{G,c}(o,\infty)$, the deletion identity becomes
\[
 \E\left[\sum_{e\in\Ray_o(F)}\Phi(F-e,e)\right]
 =r_G\E^{o,\mathrm w}\left[\sum_{e\in\partial_GA_o}c(e)\Phi(F^o,e)\right].
\]
The same one-edge inclusion coupling holds. These extensions follow from the same deletion bijection, since restoring $e$ multiplies the forest weight by $c(e)$.
\end{remark}

\subsection{The trace on the ray}\label{sec:rec}
Conditioning also on the ray identifies the walk traced onto it.
\begin{lemma}[Attachments and trace conductances]\label{lem:ray-trace}
Put $J_T(x)=j$ for $x\in \mathcal B_j$, $\tau_\gamma=\inf\{n\ge0:X_n\in\gamma\}$ and $\tau_\gamma^+=\inf\{n\ge1:X_n\in\gamma\}$. Every vertex of $H$ hits $\gamma$ almost surely, and
\begin{equation}\label{hitbush}
 \PP(J_T(x)=j\mid H,\gamma)=P_x^H(X_{\tau_\gamma}=r_j).
\end{equation}
We omit the holding mass at each ray vertex, which leaves recurrence unchanged. For $i\ne j$, the trace conductances satisfy
\begin{equation}\label{eq:trace}
 a(i,j)=\deg_H(r_i)P_{r_i}^H(X_{\tau_\gamma^+}=r_j)
 =\E[N_{ij}\mid H,\gamma],\quad
 N_{ij}=\#\{x\sim_G r_i:x\in \mathcal B_j\}.
\end{equation}
They are symmetric, have row sums at most $2d$, and satisfy $a(i,i+1)\ge1$. The trace and $H$ have the same recurrence type. Sampling instead on the lattice gives
\begin{equation}\label{eq:lattice-trace}
 b_\gamma(i,j):=\sum_{x\sim_G r_i}P_x^G(X_{\tau_\gamma}=r_j,\tau_\gamma<\infty)
 =\E[N_{ij}\mid\gamma],\quad \E[a(i,j)]=\E[b_\gamma(i,j)].
\end{equation}
\end{lemma}
\begin{proof}
Fix $H$ in the full-measure set where \eqref{eq:resampling} holds and $\mathbf Q_H$ is supported on one-ended spanning trees. If $H$ is transient, run Wilson's algorithm on $H$ rooted at infinity, starting at $o$. The first path is $\gamma$. A second walk avoiding $\gamma$ forever would create a second tree, contradicting \eqref{eq:resampling}; hence every vertex hits $\gamma$, and Wilson sampling rooted at $\gamma$ gives \eqref{hitbush}.

If $H$ is recurrent, let $V_n\uparrow V(H)$ be finite connected sets, let $H_n^{\rm w}$ be $H$ with $V(H)\setminus V_n$ identified to a vertex $z_n$, and let $T_n$ be the uniform spanning tree of $H_n^{\rm w}$ and $\gamma_n$ its path from $o$ to $z_n$. Wilson's algorithm rooted at $z_n$, starting at $o$, shows that for $x\in V_n$,
\begin{equation}\label{eq:finite-hitbush}
 \PP(J_{T_n}(x)=j\mid\gamma_n)=h_n(x,j):=P_x^{H_n^{\rm w}}(X_{\tau_{\gamma_n}}=r^n_j),
\end{equation}
where $r^n_j$ is the $j$th vertex of $\gamma_n$. Since $T_n\to T$ in law locally \cite{BLPS}, we may couple so that the edge sets agree on each finite set eventually, almost surely. We first show that the rays converge locally. Fix a finite ball $B$ of $H$. The tree $T$ is one-ended, so every edge $e\in T$ with an endpoint in $B$ has a finite side $K_e$ in $T-e$, and $e\in\gamma$ exactly when $o\in K_e$. Choose a finite ball $B'\supset B$ containing these finitely many sets $K_e$ together with their edge boundaries. Once $T_n$ and $T$ agree on the edges of $B'$ and $B'\subset V_n$, the set $K_e$ is still a side of $T_n-e$ not containing $z_n$, because every other edge of $\partial K_e$ is absent from $T_n$. Hence $e\in\gamma_n$ exactly when $o\in K_e$, with the same orientation. Thus, almost surely, $\gamma_n$ and $\gamma$ eventually contain the same edges within $B$. Since $\gamma$ is a simple infinite path, its initial segment up to its last visit to $B$ is finite, and taking $B$ larger shows that every initial segment of $\gamma_n$ eventually agrees with that of $\gamma$, including vertex indices.

Fix $x$, $j$ and $\varepsilon>0$. By recurrence, there is a finite ball $B$ such that the walk from $x$ hits $o$ before leaving $B$ with probability at least $1-\varepsilon$. On this event, the first hit on $\gamma$ or on $\gamma_n$ occurs inside $B$, the walks on $H$ and $H_n^{\rm w}$ agree before leaving $V_n\supset B$, and the hit point and its index are determined by the initial segments of the rays up to their last visits to $B$. Hence $|h_n(x,j)-P_x^H(X_{\tau_\gamma}=r_j)|\le2\varepsilon$ eventually. Similarly $J_{T_n}(x)=J_T(x)$ eventually, since the tree path from $x$ to $\gamma$ is finite. For a bounded function $\Psi$ of the edges of $\gamma$ in a finite set, \eqref{eq:finite-hitbush} gives $\E[\Psi(\gamma_n)\ind_{\{J_{T_n}(x)=j\}}]=\E[\Psi(\gamma_n)h_n(x,j)]$. Letting $n\to\infty$ by bounded convergence, then $\varepsilon\downarrow0$, gives the same identity for $\gamma$ and $P_x^H(X_{\tau_\gamma}=r_j)$. Such functions generate the $\sigma$-field of $\gamma$, which proves \eqref{hitbush}.

Decomposing the first step from $r_i$ and using \eqref{hitbush} gives \eqref{eq:trace}. Reversal of each finite excursion gives symmetry. The degree bounds the row sum, and the one-step move to $r_{i+1}$ contributes one. As the trace returns to the ray almost surely, it records every visit of the original walk to $o$. After suppressing holds its transition probabilities are $a(i,j)/\sum_{k\ne i}a(i,k)$, proving the recurrence assertion. Wilson sampling on $G$ and conditional expectation give \eqref{eq:lattice-trace}.
\end{proof}

A loop count controls the trace conductances.
\begin{lemma}[Loop bound]\label{sec:loop}\label{lem:unmarked}
For $d>4$ and $i<j$,
\begin{equation}
 \E a(i,j)\le2dZ_1\sum_{t\ge j-i+1}p_t^G(o,o)\le C(j-i)^{1-d/2},\label{eq:loopbound}
\end{equation}
where $Z_1=\sum_{t\ge0}(t+1)p_t^G(o,o)<\infty$.
\end{lemma}
\begin{proof}
Generate $\gamma=\LE(S)$ from a simple lattice walk $S$ started at $o$, and let $t_i$ be its last visit time to $r_i$. For $q=j-i$, consider pairs $(\alpha,\beta)$ in which $\alpha$ runs from $z$ to $w$ with $|\LE(\alpha)|=q$, and $\beta$ runs from $w$ to a neighbor $x$ of $z$, avoiding $\LE(\alpha)$ after time zero. Allow $|\beta|=0$. After $t_i$, the walk avoids $\gamma[0,i]$, so $\LE(S[t_i,t_j])=\gamma[i,j]$. Thus $z=r_i$, $\alpha=S[t_i,t_j]$ and the reversed first-hit path $\beta$ account for every summand of $b_\gamma(i,j)$. Enlarge this to all deterministic times $s$ with $|\LE(S[0,s])|=i$. The Markov property at these deterministic times, followed by Tonelli, gives
\begin{equation}\label{eq:prefix}
 \E b_\gamma(i,j)\le \E K_i\sum_{(\alpha,\beta)}\omega(\alpha)\omega(\beta),\qquad
 K_i=\#\{s:|\LE(S[0,s])|=i\},\quad \omega(\eta)=(2d)^{-|\eta|},
\end{equation}
with the pair sum rooted at $o$. By translation invariance that sum has the same value for each prefix endpoint $S_s$, so it factors out of the expectation. By \cite[Proposition~2.8]{LPS}, $\E K_i\le Z_1$. Map each pair to the rooted loop $\xi=\alpha\oplus\beta\oplus(x,o)$. After removing its last edge, vertex $q$ of its loop erasure is $w$, and the loop's last visit to $w$ recovers the split: the later part of $\beta$ avoids every retained vertex of $\alpha$. Thus the map is injective, including when $\beta$ is empty. The total loop weight of length $t$ is $p_t^G(o,o)$, and the closing edge costs a factor $2d$. This proves \eqref{eq:loopbound}.
\end{proof}

\subsection{Localized voltages}
Harmonic extension gives low-energy voltages; controlling their level sets produces finite-support test functions.
\begin{lemma}[Localized voltages]\label{sec:cutoffs}\label{sec:spatial}
For a finitely supported profile $f:\mathbb N_0\to[0,1]$, write $u_f(x)=E_x^H[f(J_T(X_{\tau_\gamma}))]$ for its harmonic extension, and abbreviate it to $u$ when $f$ is specified. Then
\begin{equation}\label{eq:trace-form}
 \mathcal E_H(u)=\sum_{i<j}a(i,j)(f(i)-f(j))^2.
\end{equation}
For an integer $N\ge4$, let $f_N$ equal one up to $\lfloor N/2\rfloor$, decrease linearly to zero at $N$, and vanish thereafter. Its extension $u_N=u_{f_N}$ satisfies
\begin{equation}\label{energies}
 \E\mathcal E_H(u_N)\le Ck_d(N),\quad
 k_d(N)=\begin{cases}N^{1/2}&d=5,\\1&d=6,\\N^{-1/2}&d=7,\\N^{-1}\log N&d=8,\\N^{-1}&d\ge9.\end{cases}
\end{equation}

In dimension six, the extension of
$f_N^{\rm log}(i)=(1-\log(i+1)/\log(N+1))_+$ satisfies
\begin{equation}\label{eq:d6energy}
 \E[\mathcal E_H(u_{f_N^{\rm log}})]\le C/\log(N+1).
\end{equation}
For $d\ge7$, the extension of $g_N(i)=(1-i/N)_+$ satisfies
\begin{equation}\label{eq:d7energy}
 \E[\mathcal E_H(u_{g_N})]\le CN^{-1/2}.
\end{equation}
For the plateau or logarithmic profile and an $H$-measurable vertex set $D\subset H$ containing $o$, put
$q_D=\PP(\bigcup_{j<N}\mathcal B_j\not\subset D\mid H,\gamma)$. Then $\sup_{x\notin D}u(x)\le q_D$, and
\begin{equation}\label{supportfail}
 \PP(q_D>1/4)\le4\E q_D\le C(N/R^2+e^{-cR^2/N})\qquad(D=H\cap Q_R).
\end{equation}
On $q_D\le1/4$, $(2u-1)_+$ has value one at $o$, vanishes outside $D$, and has energy at most $4\mathcal E_H(u)$. Also
\begin{equation}\label{plateaumass}
 \E\sum_{x\in H\cap Q_s}(1-u_N(x))\le Cs^dN^{(4-d)/2}.
\end{equation}
\end{lemma}
\begin{proof}
Choose finite sets $V_m\uparrow H$ containing the support of $f$ on the ray. Let $u_m$ equal $f$ on $\gamma$, vanish off $V_m$, and be harmonic at the remaining vertices. The hitting representation and $\tau_\gamma<\infty$ give $u_m(x)\to u(x)$ by bounded convergence. For $n\ge m$, the function $u_m-u_n$ vanishes on $\gamma$ and off $V_n$. Harmonicity of $u_n$ therefore gives
\[
 \mathcal E_H(u_m-u_n)=\mathcal E_H(u_m)-\mathcal E_H(u_n).
\]
Thus the gradients converge in $\ell^2(E(H))$, with pointwise limit $\nabla u$. Summing by parts for $u_m$ leaves only the finite sum $\sum_{i:f(i)\ne0} f(i)\deg_H(r_i)(I-P_H)u_m(r_i)$. Pass to the limit and use the first-step decomposition at $r_i$. The result is $\sum_i f(i)\sum_{j\ne i}a(i,j)(f(i)-f(j))$, which equals \eqref{eq:trace-form} by symmetry.

The increment bound $|f_N(i+h)-f_N(i)|\le\min(1,2h/N)$ and \eqref{eq:loopbound} give
$\E\mathcal E_H(u_N)\le CN\sum_{h\ge1}h^{1-d/2}\min(1,2h/N)^2$, which is \eqref{energies}. The same calculation gives \eqref{eq:d7energy}. For the logarithmic profile in dimension six, put $x=i+1$ and split the sum in $h$ at $x$ and then into dyadic ranges. This gives
\[
 \sum_{h\ge1}h^{-2}\log^2(1+h/x)\le C/x,\qquad
 \E[\mathcal E_H(u_{f_N^{\rm log}})]
 \le\frac{C}{\log^2(N+1)}\sum_{x\le N}x^{-1},
\]
proving \eqref{eq:d6energy}.

Conditional on $\gamma$, Lemmas~\ref{lem:past-estimates} and~\ref{lem:bush-domination} bound the intrinsic and extrinsic radius tails of each bush by $C/s$ and $C/s^2$. Lemma~\ref{lem:ray-length} keeps the first $N$ ray vertices in $Q_{R/2}$ except with probability $Ce^{-cR^2/N}$. A union bound and Markov's inequality give \eqref{supportfail}. Equation~\eqref{hitbush} gives the exterior supremum and clipping claim. Finally, attachment after $N/2$ requires a walk from $x$ to meet the independent walk from $o$ after time $N/2$. Its probability is at most $\sum_{j\ge N/2,\,k\ge0}p_{j+k}^G(o,x)\le CN^{(4-d)/2}$; summing over $Q_s$ proves \eqref{plateaumass}.
\end{proof}

\subsection{Counting components under exterior conditioning}
The packing heuristic of Section~\ref{sec:bkps} holds for the ordinary forest, and it survives conditioning on distant edges. The input is the local volume bound of Barlow and J\'arai in a box with wired boundary.
\begin{lemma}[Local volume in a wired box {\cite[Theorem~5.4]{BJ}}]\label{lem:wired-volume}
Let $d\ge5$, $h\ge1$ and $Q_{4h}(x)\subseteq D\subset\Z^d$ with $D$ finite. Let $U_D$ be the uniform spanning forest of $D$ wired at $\Z^d\setminus D$, that is, the uniform spanning tree of the graph obtained from $\Z^d$ by identifying $\Z^d\setminus D$ to one vertex and deleting the loops, with that vertex removed. If $U_D(x)$ is the component of $x$, then for $\lambda>0$,
\begin{equation}\label{eq:wired-volume}
 \PP\big(|Q_h(x)\cap U_D(x)|\le\lambda h^4\big)\le C\exp(-c\lambda^{-1/3}).
\end{equation}
\end{lemma}
This is \cite[Theorem~5.4]{BJ}, translated to $x$; for $\lambda>1$ increase $C$. The tree paths of $U_D$ between vertices of $D$ use only edges with both endpoints in $D$.

\begin{figure}[ht]
\centering
\begin{tikzpicture}[scale=0.4,font=\small,
  inf/.style={treenavy,line width=0.8pt,-{Stealth[length=4pt]}}]
 \fill[shadeA] (-4,-4) rectangle (4,4);
 \draw[rayteal,line width=0.6pt] (-4,-4) rectangle (4,4);
 \draw[rayteal,dashed,line width=0.5pt] (-2,-2) rectangle (2,2);
 \draw[treenavy!70,line width=0.6pt] (-8.5,-8.5) rectangle (8.5,8.5);
 \node[rayteal,anchor=south west,inner sep=1pt] at (-4,4) {$Q_{2s}$};
 \node[rayteal,anchor=south west,inner sep=1pt] at (-2,2) {$Q_{s}$};
 \node[treenavy,anchor=south] at (0,9.2) {boundary layer of $Q_{5s+5}$};
 \fill[extraorange!25] plot[smooth cycle,tension=0.8] coordinates {(0,0) (1.2,0.8) (2.9,1.5) (3.3,0.4) (1.6,-0.5) (0.6,-1.4) (-0.8,-0.9)};
 \draw[extraorange,line width=0.7pt] plot[smooth cycle,tension=0.8] coordinates {(0,0) (1.2,0.8) (2.9,1.5) (3.3,0.4) (1.6,-0.5) (0.6,-1.4) (-0.8,-0.9)};
 \node[circle,fill=extraorange,inner sep=1.4pt,label={[inner sep=1pt]left:$o$}] at (0.2,-0.3) {};
 \node[extraorange!80!black,anchor=south] at (2.9,1.6) {$A_o$};
 \draw[inf] (3.5,1.6) .. controls (6,3) and (8,6) .. (10.5,7);
 \draw[inf] (-0.9,-1.1) .. controls (-3,-5) and (-6,-8) .. (-7,-10.5);
 \draw[inf] (1.0,-1.6) .. controls (3,-4) and (6,-8) .. (10.5,-8);
 \draw[inf] (-0.9,0.4) .. controls (-4,2) and (-8,1) .. (-10.5,3);
 \foreach \y in {-8,-7,...,8}{\draw[latgray!60!black,line width=0.3pt] (8.05,\y-0.45) rectangle (8.95,\y+0.45); \draw[latgray!60!black,line width=0.3pt] (-8.95,\y-0.45) rectangle (-8.05,\y+0.45);}
 \foreach \x in {-7,-6,...,7}{\draw[latgray!60!black,line width=0.3pt] (\x-0.45,8.05) rectangle (\x+0.45,8.95); \draw[latgray!60!black,line width=0.3pt] (\x-0.45,-8.95) rectangle (\x+0.45,-8.05);}
 \node[anchor=west,align=left] at (9.6,-2) {cubes of\\ radius $\lceil s/10\rceil$};
\end{tikzpicture}
\caption{Lemma~\ref{lem:packing}(iii). If the wired origin tree $A_o$ reaches distance $s$, every infinite tree meeting $Q_{2s+1}$ crosses the boundary layer of $Q_{5s+5}$. This layer is covered by boundedly many cubes of radius $\lceil s/10\rceil$ whose eightfold enlargements miss $Q_s$, while the event $\{\rho\ge s\}$ depends only on edges inside $Q_s$. Not to scale.}\label{fig:decoupling}
\end{figure}
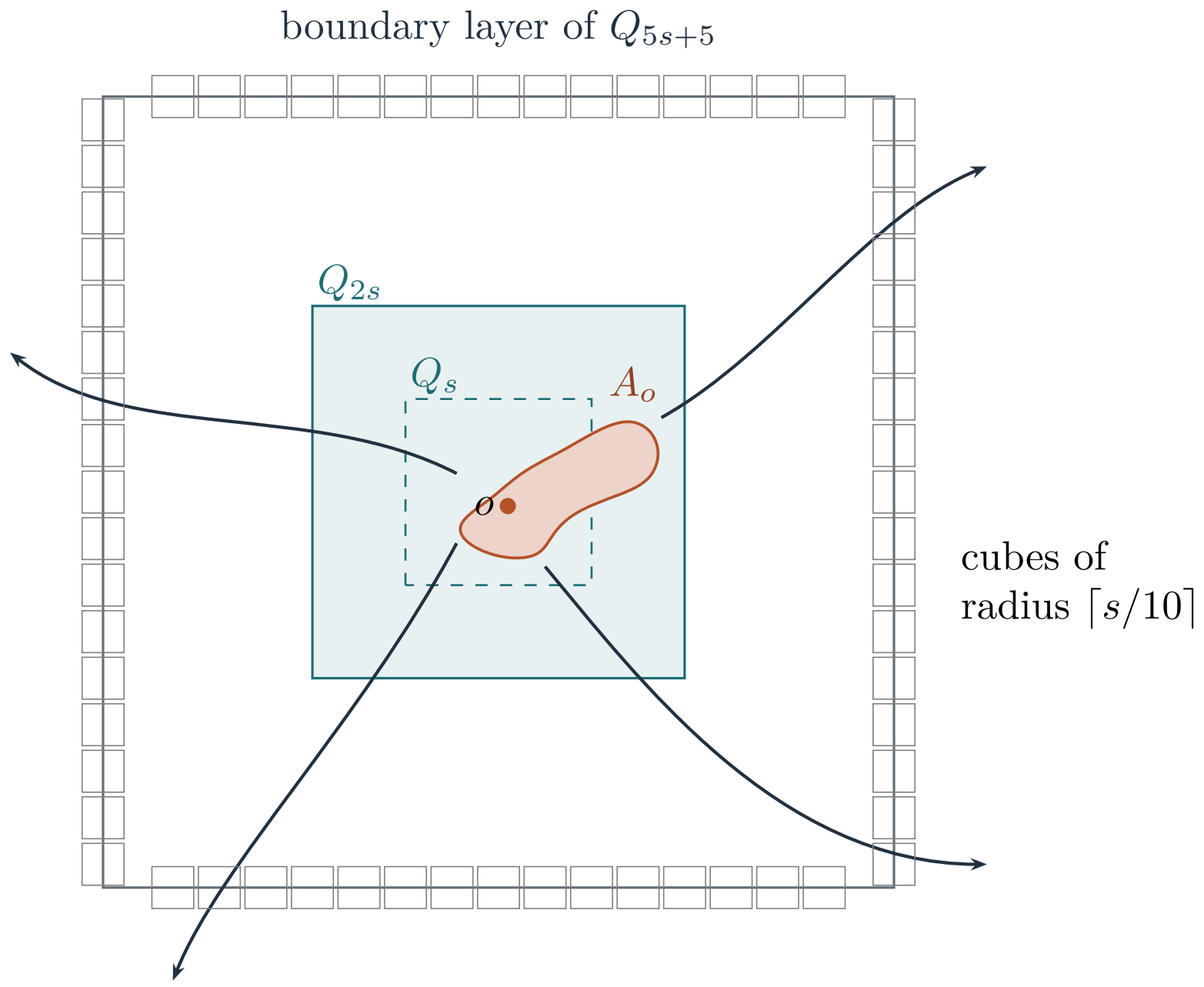

\begin{lemma}[Conditional packing]\label{lem:packing}
Let $K(S)$ and $K_\infty(S)$ count the components and the infinite components of $F^o$ meeting $S$, and $K_F(S)$ the components of $F$ meeting $S$. Put $\rho=\max_{x\in A_o}\|x\|_\infty$. For integer $s\ge1$ and $z\in\Z^d$, the following bounds hold.
\begin{enumerate}[label=\textup{(\roman*)},leftmargin=*]
\item \begin{equation}\label{eq:ordinary-packing}
 \E[K_F(Q_s(z))]\le Cs^{d-4};
\end{equation}
\item if $Q_{8s}(z)$ misses $o$ and $E$ is a cylinder event determined by edges with both endpoints outside $Q_{8s}(z)$, then
\begin{equation}\label{eq:packingclean}
 \E^{o,\mathrm w}[K(Q_s(z))\ind_E]\le Cs^{d-4}\PP^{o,\mathrm w}(E);
\end{equation}
\item and
\begin{equation}\label{eq:conditionalD}
 \E^{o,\mathrm w}[K_\infty(Q_{2s+1})\ind_{\{\rho\ge s\}}]\le Cs^{d-4}\PP^{o,\mathrm w}(\rho\ge s)\le Cs^{d-6}.
\end{equation}
\end{enumerate}
Every component of $F^o$ other than $A_o$ is infinite.
\end{lemma}
\begin{proof}
\emph{Local volume under exterior conditioning.} For $x$ with $Q_{4h+1}(x)\subset Q_{8s}(z)$, let $M_h(x)$ count the vertices of $Q_h(x)$ joined to $x$ by forest edges with both endpoints in $Q_{4h}(x)$. We claim that
\begin{equation}\label{eq:local-packing-tail}
 \PP^{o,\mathrm w}(E,\ M_h(x)\le\lambda h^4)\le C\PP^{o,\mathrm w}(E)\exp(-c\lambda^{-1/3}),
\end{equation}
and the same bound for $F$ with $E$ trivial. Work in a finite wired box containing $Q_{8s}(z)$, with $o$ identified with the sink for $\PP^{o,\mathrm w}$; both sides are limits of their finite-volume versions, since $E$ is a cylinder event and $M_h(x)$ is a function of finitely many edges. Condition on the states of all edges with both endpoints outside $Q_{8s}(z)$. By the spatial Markov property, the remaining edges form the uniform spanning tree of the graph in which the present conditioned edges are contracted and the absent ones deleted. Identifying, in addition, all vertices outside $Q_{4h}(x)$ with the sink adds roots, so by Lemma~\ref{lem:wilson-facts} it decreases the remaining forest in stochastic order. The edges of the resulting forest inside $Q_{4h}(x)$ have the law of $U_D$ with $D=Q_{4h}(x)$. Since $M_h(x)$ is increasing in the edge set, Lemma~\ref{lem:wired-volume} bounds the conditional probability of $\{M_h(x)\le\lambda h^4\}$ by $C\exp(-c\lambda^{-1/3})$. Integrating over the conditioning on $E$ proves \eqref{eq:local-packing-tail}.

\emph{Counting witnesses.} Put $W_j=Q_{s+s/(j+1)}(z)$, $t_j=2^{-j}s^4$, $h_j=\lfloor s/[2j(j+1)]\rfloor$ for $j\ge1$, and $J=\lceil4\log_2s\rceil+1$. Assign each tree $\mathcal T$ meeting $Q_s(z)$ to the least $j\in\{0,\ldots,J\}$ with $|\mathcal T\cap W_j|\ge t_j$; since $t_J<1$ and $Q_s(z)\subset W_J$, such $j$ exists. The trees with index zero have disjoint subsets of $W_0=Q_{2s}(z)$ of size at least $s^4$, so there are at most $Cs^{d-4}$ of them. A tree with index $j\ge1$ satisfies $|\mathcal T\cap W_{j-1}|<2t_j$. Since $s/j-s/(j+1)\ge h_j$, each of its at least $t_j$ vertices $x\in\mathcal T\cap W_j$ has $Q_{h_j}(x)\subset W_{j-1}$, hence $M_{h_j}(x)\le|\mathcal T\cap Q_{h_j}(x)|<2t_j$. Moreover $Q_{4h_j+1}(x)\subset Q_{8s}(z)$. For large $s$ and $j\le J$ we have $h_j\ge s/[4j(j+1)]\ge1$, so $h_j^4/(2t_j)\ge c2^j/(j+1)^8$. By \eqref{eq:local-packing-tail}, the expected number of trees with index $j$, on $E$, is at most
\[
 t_j^{-1}\sum_{x\in Q_{2s}(z)}\PP^{o,\mathrm w}\big(E,\ M_{h_j}(x)<2t_j\big)
 \le Cs^{d-4}\PP^{o,\mathrm w}(E)\,2^j\exp\{-c[2^j/(j+1)^8]^{1/3}\}.
\]
Summation over $j$ proves \eqref{eq:packingclean}; bounded $s$ are absorbed into the constant. The same argument for $F$ without $E$ proves (i).

\emph{Infinite neighbors.} The event $\{\rho\ge s\}$ is determined by edges with both endpoints in $Q_s$: it occurs exactly when $o$ is joined to a vertex of norm $s$ by a forest path inside $Q_s$. Every infinite component meeting $Q_{2s+1}$ leaves $Q_{5s+5}$, so it meets the boundary layer $Q_{5s+5}\setminus Q_{5s+4}$. For $s\ge10$, cover this layer by at most $C_d$ cubes $Q_{\lceil s/10\rceil}(z_i)$ centered in it (Figure~\ref{fig:decoupling}); their eightfold enlargements miss $Q_s$, and in particular miss $o$. Apply \eqref{eq:packingclean} to each cube with $E=\{\rho\ge s\}$, and use $\PP^{o,\mathrm w}(\rho\ge s)\le Cs^{-2}$ from \eqref{eq:Htails}. This proves \eqref{eq:conditionalD}; bounded $s$ are absorbed into the constant.

Finally, every component other than $A_o$ is infinite: a finite component avoiding $o$ would give a cylinder event absent from every sufficiently large finite approximation.
\end{proof}

\section{Recurrence and transience}\label{sec:five}
Deletion and component counting prove transience in dimension five; the trace-energy bounds prove recurrence in higher dimensions. We then record a cut formula for box resistance and identify its scale in dimensions six and seven.
\subsection{Recurrence classification}\label{sec:recurrence}\label{subsec:contact}
We use deletion for transience and the trace capacity for recurrence.
\begin{proof}[Proof of Theorem~\ref{thm:main}]
On $s\le\rho<2s$, every component adjacent to $A_o$ meets $Q_{2s+1}$, and it is infinite by Lemma~\ref{lem:packing}. Thus \eqref{eq:conditionalD} gives $\E^{o,\mathrm w}[D_o\ind_{\{s\le\rho<2s\}}]\le Cs^{d-6}$. In dimension five, summing over dyadic $s\ge R$ yields
\begin{equation}\label{eq:five-contact-tail}
 \E^{o,\mathrm w}[D_o\ind_{\{\rho\ge R\}}]\le C/R,\qquad \E^{o,\mathrm w}D_o<\infty.
\end{equation}
The contribution at $\rho=0$ is at most $2d$. The deletion identity now proves finite mean escape resistance and almost-sure transience.

For recurrence, use the finitely supported ray profile $f_N^{\rm log}$ in dimension six and $g_N$ in dimensions $d\ge7$. Both equal one at index zero. Write $\operatorname{Cap}_a(0)$ for the infimum of $\sum_{i<j}a(i,j)(f(i)-f(j))^2$ over finitely supported $f:\mathbb N_0\to\mathbb R$ with $f(0)=1$. By \eqref{eq:trace-form}, \eqref{eq:d6energy} and \eqref{eq:d7energy}, this capacity satisfies
\[
 \E[\operatorname{Cap}_a(0)]\le
 \begin{cases}C/\log(N+1)&d=6,\\ CN^{-1/2}&d\ge7.\end{cases}
\]
Letting $N\to\infty$ gives $\operatorname{Cap}_a(0)=0$ almost surely. After its holding steps are suppressed, the trace chain jumps from $i$ to $j\ne i$ with probability $a(i,j)/\sum_{k\ne i}a(i,k)$. These denominators are finite and positive. The capacity criterion \cite[Theorem~2.3]{LP} gives recurrence of this chain, and Lemma~\ref{lem:ray-trace} transfers it to $H$. Translation invariance and a countable intersection over roots give the assertion for all components.
\end{proof}

\subsection{A cut formula for box resistance}\label{sec:terminal}
Stopping the cut potential at an exit time gives a finite-target version of the deletion formula.
\begin{lemma}[Cut formula for an exit set]\label{lem:terminal}
In every recurrence type, with $\tau_R=\inf\{n\ge0:X_n\notin Q_R\}$,
\begin{equation}\label{terminal}
 \E\mathcal R_R=r_G\E^{o,\mathrm w}\sum_{K\sim A_o}
 P_o^{\Z^d[A_o\cup K]}(X_{\tau_R}\in K).
\end{equation}
Here the sum ranges over components $K$ of $F^o$ other than $A_o$ with $m(A_o,K)>0$. In dimensions six and seven, the box resistance satisfies
\begin{equation}\label{eq:six-mean-resistance}
 \mathcal R_R\asymp_{\PP} b_d(R),\qquad \E\mathcal R_R\asymp b_d(R),
 \qquad b_6(R)=\log(2+R),\quad b_7(R)=R.
\end{equation}
\end{lemma}
\begin{proof}
For a finite connected $A\subset H$, let $L_A$ be the event that $T[A]$ is connected and has exactly one forest exit. For $z\in H$, put $a_o^H(z)=\sum_{A\ni o,\,z\notin A}\mathbf Q_H(L_A)/b_H(A)$, summing over these sets. Every such set crosses a fixed path from $o$ to $z$; by uniform exit, its weights at any path edge sum to at most the probability that the ray uses that edge. Hence $0\le a_o^H(z)\le d_H(o,z)$. Its increments are the signed ray current, with $(P_H-I)a_o^H=\ind_o/\pi(o)$. From every vertex of the finite set $H\cap Q_R$ there is a path to its exterior, so $\tau_R$ has a geometric tail for fixed $H,R$. The potential is bounded on the finite set $H\cap Q_{R+1}$. Stopping the Poisson equation at $\tau_R$ gives
\[
 \mathcal R_R=\frac1{\pi(o)}E_o^H\left[\sum_{n<\tau_R}\ind_{\{X_n=o\}}\right]
 =E_o^H[a_o^H(X_{\tau_R})].
\]
Apply \eqref{eq:marked} to each cut and its exit probability. In the deleted configuration all $m(A_o,K)$ marks give the same induced graph and probability, each with weight $1/m(A_o,K)$. Tonelli proves \eqref{terminal}. On $\rho\ge R$, every contributing component meets $Q_{R+1}$, so \eqref{eq:conditionalD} bounds its expected contribution by $CR^{d-6}$. On $\rho<R$, sum the shell bounds $Cs^{d-6}$ used in \eqref{eq:five-contact-tail}, including the bounded $\rho=0$ contribution. This yields $\E\mathcal R_R\le Cb_d(R)$ for $d=6,7$.

For the converse, fix $\eta>0$ and set $D=H\cap Q_R$ and $N=\lfloor aR^2\rfloor$. Choose $a>0$ small enough that \eqref{supportfail} gives $\PP(q_D>1/4)<\eta/2$ for large $R$. Use the logarithmic profile in dimension six and the plateau profile in dimension seven. By \eqref{eq:d6energy} and \eqref{energies}, their extensions have mean energy at most $C_a/b_d(R)$. Markov's inequality gives $\mathcal E_H(u)\le 2C_a/(\eta b_d(R))$ outside an event of probability $\eta/2$. On the intersection of the two good events, $(2u-1)_+$ is admissible for the box capacity, so Dirichlet's principle gives
$\mathcal R_R\ge (4\mathcal E_H(u))^{-1}\ge c_\eta b_d(R)$. This proves the lower bounds in \eqref{eq:six-mean-resistance}. Its probabilistic upper bound follows by Markov's inequality from the mean upper bound.
\end{proof}

\section{Return probabilities and fluctuations}\label{sec:heat}
We prove the return estimates and fluctuations in Theorem~\ref{thm:returns}, apart from the sharp eight-dimensional bounds completed in Section~\ref{sec:d8obs}.

Section~\ref{sec:heatlower} obtains return lower bounds from localized voltages, using the Hilbert-space energy and mass-transport method of \cite{Lee23,LeeScaling}. The required model-specific input is a family with the correct energy and sufficiently many nonzero coordinates near the root. Section~\ref{sec:heatupper} proves the upper bounds from a cut representation, whose finite-network form follows from the two-component forest Green formula of \cite[Theorem~3]{KKW}. We estimate small cuts through the Laplacian and large cuts through the energy; the conditional cut counts make this decomposition sharp in dimensions five through seven, without a logarithmic loss in dimension six.

For the fluctuations in Section~\ref{sec:heat6}, the use of atypical Wilson paths has precedents in \cite{BCK,SW}. Here we tilt only a finite middle segment of the first walk, lowering the voltage energy on the induced graph while retaining a uniformly positive amount of mass near the root. Section~\ref{sec:d6} collects the remaining six-dimensional estimates.

Throughout, $\pi(x)=\deg_H(x)$ and $P=P_H$; an omitted superscript in $p_n$ refers to the walk on $H$.

\subsection{Lower bounds for return probabilities}\label{sec:heatlower}
Put $s_d=\min\{d-6,2\}$ and $w_d=4+s_d$. Translation-covariant localized voltages confine the walk to a lattice box; a volume upper bound then gives the return lower bound.

\begin{lemma}[Return lower bounds]\label{lem:returnlower}
For $d\ge5$, almost surely
\begin{equation}\label{eq:heatlower}
 \liminf_{n\to\infty}\frac{\log p_{2n}^H(o,o)}{\log n}\ge-4/w_d.
\end{equation}
For every $\eta>0$, at every sufficiently large deterministic $n$,
\begin{equation}\label{eq:typicallower}
 \PP(p_{2n}^H(o,o)\ge c_\eta a_d(n))\ge1-\eta,
 \qquad a_d(n)=\begin{cases}n^{-4/w_d}&d=5,6,7\text{ and }d\ge9,\\
 (n\log(2+n))^{-2/3}&d=8.
 \end{cases}
\end{equation}
In particular $\E p_{2n}^H(o,o)\ge c\,a_d(n)$.
\end{lemma}
\begin{proof}
We first prove the confinement bound \eqref{eq:confinement-general}, then choose scales for the two assertions.

\emph{Localized voltages and active centers.}
Fix $N\ge4$ and extend the plateau profile $f_N$ from the ray at each $z$ by \eqref{hitbush}, obtaining $u_z$ on $H_z$. Declare $z$ bad if $\sup_{H_z\setminus Q_R(z)}u_z>1/4$, and set $h_z=(u_z-1/4)_+$ at good centers and zero otherwise, extending it by zero off $H_z$. Thus $h_z$ is supported in $Q_R(z)$, and the family is translation covariant. Equations \eqref{energies} and \eqref{supportfail} give
\begin{equation}\label{eq:centerenergy}
 \PP(z\text{ bad})\le e_0:=C(N/R^2+e^{-cR^2/N}),\qquad
 \E\mathcal E_{H_z}(h_z)\le Ck_d(N).
\end{equation}
For $r<R/2$, let $Z$ count bad centers in $H\cap Q_r$ and put $B=\sum_{z\in H\cap Q_r}(1-u_z(o))$. For $p\ge2$, mass transport, H\"older, Lemma~\ref{lem:forest-volumes}, and \eqref{plateaumass} give
\[
 \begin{gathered}
 \E Z=\E[\ind_{\{o\text{ bad}\}}V_r]\le C_pr^4e_0^{1-1/p},\\
 \E B=\E\sum_{x\in H\cap Q_r}(1-u_o(x))\le Cr^dN^{(4-d)/2}.
 \end{gathered}
\]
For $0<v<1$, the volume lower tail and Markov's inequality give $V_r\ge vr^4$ and $Z+2B\le vr^4/2$, except on probability at most
\[
 \eta_0=Ce^{-cv^{-1/3}}+C_pv^{-1}
 \big[e_0^{1-1/p}+(r/\sqrt N)^{d-4}\big].
\]
Every good center with $u_z(o)\ge1/2$ has $h_z(o)\ge1/4$; the other centers number at most $Z+2B$. Hence, outside the same exceptional event,
\begin{equation}\label{eq:manycenters}
 \#\{z\in H\cap Q_r:h_z(o)\ge1/4\}\ge vr^4/2.
\end{equation}
\emph{From energy to confinement.}
The vectors $\Phi(x)=(h_z(x))_{z\in\mathbb Z^d}$ belong to $\ell^2(\mathbb Z^d)$ and have finite support, since $h_z(x)=0$ for $z\notin Q_R(x)$. If $X_m\notin Q_{2R}$, the coordinates counted in \eqref{eq:manycenters} vanish at $X_m$, so
$\|\Phi(X_m)-\Phi(o)\|^2\ge vr^4/32$.
Mass transport exchanges the root with the coordinate center. By reversibility,
\begin{equation}\label{eq:voltage-displacement}
 \begin{aligned}
 \E[\pi(o)E_o^H\|\Phi(X_m)-\Phi(o)\|^2]
 &=\E\sum_{x\in H}\pi(x)E_x^H(h_o(X_m)-h_o(x))^2\\
 &=2\E\langle h_o,(I-P_H^m)h_o\rangle_\pi
 \le2m\E\mathcal E_H(h_o).
 \end{aligned}
\end{equation}
The inequality follows from $1-t^m\le m(1-t)$ for $-1\le t\le1$. Since $\pi(o)\ge1$, the same upper bound controls the unweighted expectation. Markov's inequality, \eqref{eq:centerenergy}, and the exceptional probability $\eta_0$ therefore yield
\begin{equation}\label{eq:confinement-general}
 \E P_o^H(X_m\notin Q_{2R})\le C\eta_0+\frac{Cm k_d(N)}{vr^4}.
\end{equation}
\emph{Return probabilities and choice of scales.}
The confinement-to-return estimate \cite[Lemma~2.3]{Lee23} gives
\begin{equation}\label{eq:confinement-return}
 p_{2m}^H(o,o)\ge
 \frac{\pi(o)P_o^H(X_m\in Q_{2R})^2}
 {\sum_{x\in H\cap Q_{2R}}\pi(x)}.
\end{equation}
For the almost-sure bound, set $p=2$ and take $N=\lfloor R^{2-\theta}\rfloor$, $r=\lfloor R^{1-\theta}\rfloor$, $v=R^{-\theta/8}$, and $m=\lfloor R^{w_d-8\theta}\rfloor$, with $0<\theta<1/10$. The bound \eqref{eq:confinement-general} is $O(R^{-3\theta/8})$, summable on dyadic radii, so the quenched escape probability tends to zero almost surely. Lemma~\ref{lem:forest-volumes} and \eqref{eq:confinement-return} then give $p_{2m}^H(o,o)\ge R^{-4-o(1)}$. The times have bounded successive ratios, so monotonicity of even returns and then $\theta\downarrow0$ prove \eqref{eq:heatlower}.

For fixed confidence, keep $p=2$, choose $v$ small, then $N=\lfloor aR^2\rfloor$ with $a$ small, and finally $r=\lfloor bR\rfloor$ with $b/\sqrt a$ small. These choices make $\eta_0$ small; choosing $c_\eta$ last makes \eqref{eq:confinement-general} arbitrarily small whenever $m\le c_\eta R^{w_d}$ for $d\ne8$, or $m\le c_\eta R^6/\log R$ for $d=8$. Markov's inequality and the volume upper tail give quenched mass at least $1/2$ and volume at most $C_\eta R^4$ with the prescribed confidence. Equation~\eqref{eq:confinement-return} gives $p_{2m}\ge c_\eta R^{-4}$; the least admissible $R$ proves \eqref{eq:typicallower}. A fixed positive confidence gives the mean lower bound.
\end{proof}

\subsection{Upper bounds from finite cuts}\label{sec:heatupper}
We extend the finite-network Green representation of \cite[Theorem~3]{KKW} to the infinite one-ended network $H$. Its cut weights give return upper bounds by controlling small cuts with the Laplacian and large cuts with energy.

\begin{lemma}[Dipoles and cut energy]\label{lem:dipoles}
For each nonempty finite connected $A\subset H$, let $L_A$ mean that $T[A]$ is connected and has exactly one forest exit, whether or not $o\in A$, and put $w_A=\mathbf Q_H(L_A)/b_H(A)$. All cut sums below range over these sets. Almost surely, for $x\in H$,
\begin{gather}
 \res^H(o,x)=\sum_Aw_A(\ind_A(o)-\ind_A(x))^2,\label{dipole}\\
 \E\sum_{x\in H\cap Q_R}\res^H(o,x)=2r_G\E^{o,\mathrm w}\sum_{K\sim A_o}|K\cap Q_R|.\label{pairavg}
\end{gather}
For finitely supported $h$, put $\mathsf{L}=I-P_H$ and $\mu(y)=\pi(y)\mathsf{L}h(y)$. Then
\begin{gather}
 h(o)-h(x)=\sum_Aw_A(\ind_A(o)-\ind_A(x))\mu(A),\label{fixedcutdiff}\\
 \mathcal E_H(h)=\sum_Aw_A\mu(A)^2.\label{fixedcutenergy}
\end{gather}
\end{lemma}
\begin{proof}
For \eqref{dipole} the case $x=o$ is immediate. Take $x\ne o$ and put $s_A=\ind_A(o)-\ind_A(x)$. For each edge of a fixed $o$--$x$ path, the uniform-exit identity bounds the sum of $w_A$ over cuts crossing it by one: a tree edge has only one finite side. Thus $\sum_Aw_A|s_A|<\infty$, and
$g(y)=\sum_Aw_As_A\ind_A(y)$ converges uniformly and vanishes at infinity. Its gradient is the mean signed current along the tree path from $o$ to $x$, so its divergence is $\ind_o-\ind_x$. Test against the finitely supported truncation $(g-t)_+-(-g-t)_+$ and let $t\downarrow0$ to obtain
$\mathcal E_H(g)=g(o)-g(x)=\sum_Aw_As_A^2$.
Thus $\nabla g$ is a unit flow of energy $g(o)-g(x)$. Conversely, write $g_t=(g-t)_+-(-g-t)_+$. For any finite-energy unit flow $i$ from $o$ to $x$, summation by parts gives $\langle i,\nabla g_t\rangle=g_t(o)-g_t(x)$. As $t\downarrow0$, the gradients converge in energy by dominated convergence. Cauchy--Schwarz now gives $\sum_e i(e)^2\ge g(o)-g(x)$, which proves \eqref{dipole}. For \eqref{pairavg}, deletion turns the contribution of cuts containing $o$ but not $x$ into
$r_G\E^{o,\mathrm w}\sum_{K\sim A_o}|K\cap Q_R|$.
Interchanging roots gives the same contribution from cuts containing $x$ but not $o$, by translation invariance and symmetry of $Q_R$.

For general $h$, express the finitely supported, zero-mass charge $\mu=\pi \mathsf{L}h$ as a finite sum of dipoles. Its potential $g(y)=\sum_Aw_A\mu(A)\ind_A(y)$ has $\pi \mathsf{L}g=\mu$ and vanishes at infinity. The maximum principle gives $g=h$. Subtracting its values at $o,x$ proves \eqref{fixedcutdiff}; summing against $\mu$ proves \eqref{fixedcutenergy}.
\end{proof}

The next lemma converts voltage bounds to return bounds; its resistance average also applies in dimension eight.
\begin{lemma}[Resolvent and resistance averaging]\label{lem:resolvent-averaging}
Let $\mathsf{L}=I-P_H$ and $\delta=\ind_o/\pi(o)$. For $\lambda>0$ and $n\ge1$,
\begin{equation}\label{resolventtime}
 G_\lambda:=\langle\delta,(\mathsf{L}+\lambda I)^{-1}\delta\rangle_\pi
 =\sup_{f\ne0}\frac{f(o)^2}{\mathcal E_H(f)+\lambda\|f\|_\pi^2},
 \qquad p_{2n}^H(o,o)\le C\pi(o)G_{1/n}/n,
\end{equation}
where the supremum is over finitely supported $f$. Writing $V_r=|H\cap Q_r|$ and $Y_r=\sum_{x\in H\cap Q_r}\res^H(o,x)$, we have
\begin{equation}\label{d8av}
 p_{2n}^H(o,o)\le C\left(V_r^{-1}+\frac{Y_r}{nV_r}\right).
\end{equation}
\end{lemma}
\begin{proof}
The variational identity follows from Cauchy--Schwarz for $\mathsf{L}+\lambda I$ and density of finitely supported functions. In the expansion
$G_\lambda=\pi(o)^{-1}\sum_{k\ge0}(1+\lambda)^{-k-1}p_k^H(o,o)$,
retain even times up to $2n$ and use monotonicity of even returns to obtain the second assertion of \eqref{resolventtime}.
For finitely supported $f$, pairing its gradient with the dipole current in Lemma~\ref{lem:dipoles} gives $(f(o)-f(x))^2\le\res^H(o,x)\mathcal E_H(f)$. For any nonempty finite $B\subset H$, average
$f(o)^2\le2f(x)^2+2\res^H(o,x)\mathcal E_H(f)$
over $x\in B$. Since $\pi(x)\ge1$, the variational identity gives
\[
 G_{1/n}\le\frac{2n}{|B|}+\frac2{|B|}\sum_{x\in B}\res^H(o,x).
\]
Take $B=H\cap Q_r$ and apply \eqref{resolventtime}, using $\pi(o)\le2d$, to prove \eqref{d8av}.
\end{proof}

In dimensions five through seven, we estimate two cut sums. The second sum below selects large cuts separating their center from a vertex in its $R$-box.
\begin{lemma}[Mean cut weights]\label{lem:cut-statistics}
Let $d=5,6,7$ and $R\ge1$. With $\pi(A)=\sum_{z\in A}\pi(z)$, define for $y\in H$
\[
 M_R^*(y)=\sum_{\substack{A\ni y\\\diam_{\rm ext} A\le R}}w_A\sqrt{\pi(A)},\qquad
 \mathcal W_R(y)=\sum_{\substack{A\ni y,\ \diam_{\rm ext} A>R\\(H\cap Q_R(y))\setminus A\ne\varnothing}}w_A.
\]
Then
\begin{equation}\label{eq:cut-statistics}
 \E M_R^*(o)\le CR^{d-4},\qquad
 \E\mathcal W_R(o)\le CR^{d-6}.
\end{equation}
\end{lemma}
\begin{proof}
For dyadic $s\ge1$, let $\mathcal A_s$ consist of the sets $A\ni o$ with
$s\le\max_{z\in A}\|z\|_\infty<2s$, including radius zero when $s=1$, and put $C_s(o)=\sum_{A\in\mathcal A_s}w_A$. Deletion and \eqref{eq:conditionalD} give
\[
 \E C_s(o)\le C\E^{o,\mathrm w}
 [K_\infty(Q_{2s+1})\ind_{\{\rho\ge s\}}]+C\ind_{\{s=1\}}
 \le Cs^{d-6}.
\]
For $\mathcal W_R$, deletion leaves only neighboring components meeting $Q_R$, while $\diam_{\rm ext} A_o>R$ implies $\rho>R/2$. Thus the same estimate, at scale $R/2$, gives
\[
 \E\mathcal W_R(o)\le C\E^{o,\mathrm w}
 [K_\infty(Q_R)\ind_{\{\rho>R/2\}}]\le CR^{d-6},
\]
with bounded radii absorbed into the constant.

On the shell $\mathcal A_s$, deletion also gives
\[
 \E\sum_{A\in\mathcal A_s}w_A|A|
 \le C\E^{o,\mathrm w}
 [|A_o\cap Q_{2s}|K_\infty(Q_{2s+2})].
\]
To bound this expectation, fix $z\in Q_{2s}$, condition on $z\in A_o$, and run Wilson first from $z$, revealing its loop-erased path to $o$. For each revealed path, the unforced edges have the forest law rooted there and at infinity; Lemma~\ref{lem:wilson-facts} couples them below an ordinary WUSF with its unconditioned marginal law. Every nonroot component is infinite. Two such components inside the same ordinary tree would give two disjoint rays in that tree, whereas rays in a one-ended tree eventually coalesce. The injection into ordinary components meeting $Q_{2s+2}$ and \eqref{eq:ordinary-packing} therefore bound the conditional expected count by $Cs^{d-4}$. Now sum
$\PP^{o,\mathrm w}(z\in A_o)=g_G(z,o)/g_G(o,o)$
over $z\in Q_{2s}$; the Green estimate in Lemma~\ref{lem:lattice-walk} makes this sum at most $Cs^2$. Hence
$\E\sum_{A\in\mathcal A_s}w_A|A|\le Cs^{d-2}$.
Since $\pi(A)\le2d|A|$, Cauchy--Schwarz for the joint environment--cut measure yields
\[
 \E\sum_{A\in\mathcal A_s}w_A\sqrt{\pi(A)}
 \le\left(\E C_s(o)\right)^{1/2}
 \left(\E\sum_{A\in\mathcal A_s}w_A\pi(A)\right)^{1/2}
 \le Cs^{d-4}.
\]
Every cut in $M_R^*(o)$ lies in a shell with $s\le R$; summing over dyadic $s$ proves the remaining assertion of \eqref{eq:cut-statistics}.
\end{proof}

Combining the cut means with the volume lower tail gives a return tail whose integrability will yield the moment estimates.
\begin{lemma}[Upper return bounds]\label{lem:returnupper}
For $d=5,6,7$,
\begin{equation}\label{eq:unified-return-tail}
 \PP\big(n^{4/(d-2)}p_{2n}>Ct\big)
 \le C[\log(2+t)]^{3(d-4)/4}t^{-(d-2)/4}
 \qquad(n,t\ge1).
\end{equation}
Consequently, for every $\epsilon>0$,
\begin{equation}\label{fixedheatupper}
 \PP(p_{2n}^H(o,o)\le C_{d,\epsilon}n^{-4/(d-2)})\ge1-\epsilon,
\end{equation}
and $p_{2n}^H(o,o)\le n^{-4/(d-2)+o(1)}$ almost surely. For $d\ge8$, almost surely $p_{2n}^H(o,o)\le n^{-2/3+o(1)}$, and $\PP(p_{2n}^H(o,o)\le C_\epsilon n^{-2/3})\ge1-\epsilon$ and $\E p_{2n}^H(o,o)\le Cn^{-2/3}$; these upper bounds are sharp for $d\ge9$, and are improved by a logarithm in Section~\ref{sec:contact-sharp} for $d=8$.
\end{lemma}
\begin{proof}
We first control returns when enough nearby vertices have bounded cut sums, then choose scales. For $d\ge8$ we average over a tree ball.

\emph{Pointwise control and smoothing.}
Fix $R,n\ge1$ and $U,V>0$. Call $y\in H$ good if $M_R^*(y)\le U$ and $\mathcal W_R(y)\le V$. For good $o,x\in H\cap Q_R$, split \eqref{fixedcutdiff} at diameter $R$. The small cuts contribute at most
$(M_R^*(o)+M_R^*(x))\|\mathsf{L}h\|_\pi$,
because $|\mu(A)|\le\sqrt{\pi(A)}\|\mathsf{L}h\|_\pi$. Every large cut separating $o,x$ is counted by either $\mathcal W_R(o)$ or $\mathcal W_R(x)$. Cauchy--Schwarz and \eqref{fixedcutenergy} therefore give
\[
 |h(o)-h(x)|\le C\{U\|\mathsf{L}h\|_\pi+\sqrt{V\mathcal E_H(h)}\}.
\]
If $o$ and at least half of $H\cap Q_R$ are good, averaging over those vertices yields
\[
 h(o)^2\le C\{V_R^{-1}\|h\|_\pi^2
       +U^2\|\mathsf{L}h\|_\pi^2+V\mathcal E_H(h)\}.
\]
Let $\widetilde P=(I+P_H)/2$, with $n$-step kernel $\widetilde p_n(x,y)$, and take $h=\widetilde P^{\,n}f$ for finitely supported $f$ with $\|f\|_\pi=1$. Spectral calculus gives
$\|h\|_\pi\le1$, $\|\mathsf{L}h\|_\pi\le C/n$ and $\mathcal E_H(h)\le C/n$.
The squared evaluation norm of $\widetilde P^{\,n}$ at $o$ is $\widetilde p_{2n}(o,o)/\pi(o)$, so
$\widetilde p_{2n}(o,o)\le C(V_R^{-1}+U^2/n^2+V/n)$.
In the binomial expansion of $\widetilde p_{2n}(o,o)$, even indices have total weight $1/2$, and their returns are all at least $p_{2n}$. Thus $p_{2n}\le2\widetilde p_{2n}(o,o)$.

\emph{The fraction of good vertices.}
For any translation-covariant bad-site event $E$, send mass $\ind_{E(x)}/V_R(u)$ from $u$ to same-component $x\in Q_R(u)$, where $V_R(u)$ is the component volume in $Q_R(u)$. Partition the senders to each recipient into boundedly many diameter-$R$ cells. Each sender sees its entire cell, whose reciprocal-volume weights therefore sum to at most one. Mass transport gives
\[
 \E\left[\frac{\#\{x\in H\cap Q_R:E(x)\}}{V_R}\right]
 \le C_d\PP(E(o)).
\]
Markov and \eqref{eq:cut-statistics} control a bad origin and a bad fraction exceeding $1/2$. With the volume lower tail, this gives, for $R,n\ge1$, $0<b<1$ and $U,V>0$,
\begin{equation}\label{eq:cut-return-tail}
 \PP\!\left(p_{2n}>C\{(bR^4)^{-1}+U^2/n^2+V/n\}\right)
 \le Ce^{-cb^{-1/3}}+C\left(\frac{R^{d-4}}U+\frac{R^{d-6}}V\right).
\end{equation}

\emph{Choice of scales.}
Put $F_n=n^{4/(d-2)}$ and $\ell=\log(2+t)$. For $1\le t\le F_n$, choose dyadic
$R\asymp(F_n/t)^{1/4}\ell^{3/4}$ and set
\[
 b=c_0\ell^{-3},\qquad U=n\sqrt{t/F_n},\qquad V=nt/F_n.
\]
The return threshold in \eqref{eq:cut-return-tail} is $Ct/F_n$. Choose $c_0$ small so that the volume error is at most $Ct^{-2}$. Since $F_n^{(d-2)/4}=n$, the two other errors are bounded respectively by
$C\ell^{3(d-4)/4}t^{-(d-2)/4}$ and
$C\ell^{3(d-6)/4}t^{-(d-2)/4}$.
This proves \eqref{eq:unified-return-tail}; for $t>F_n$ its event is empty after increasing $C$ if necessary. Taking $t$ large proves \eqref{fixedheatupper}. Taking $n=2^j$, $t=(j+1)^4$, and applying Borel--Cantelli and monotonicity proves the almost-sure upper exponents.

For $d\ge8$, take $B=B_T(o,m)$ in the averaging argument of Lemma~\ref{lem:resolvent-averaging}. Sending one unit of current along the tree path gives $\res^H(o,x)\le m$ for every $x\in B$. Thus $p_{2n}\le C(|B_T(o,m)|^{-1}+m/n)$, and with $m=\lfloor n^{1/3}\rfloor$ this gives
\[
 p_{2n}\le Cn^{-2/3}\left(1+\frac{m^2}{|B_T(o,m)|}\right).
\]
The volume tails and negative moments in Lemma~\ref{lem:intrinsic-volume} give all the asserted upper bounds.
\end{proof}

Tail integration gives the five-dimensional fractional moments and seven-dimensional mean; the return lower bound makes both sharp.
\begin{lemma}[Fractional moments in dimension five]\label{lem:momentsfive}
In $d=5$, for every $0<q<3/4$,
\[
 \E[(p_{2n}^H(o,o))^q]\asymp_q n^{-4q/3}.
\]
\end{lemma}
\begin{proof}
Integrating \eqref{eq:unified-return-tail} against $qt^{q-1}\,dt$ gives
\[
 \sup_{n\ge1}\E[(n^{4/3}p_{2n})^q]
 \le C_q\left(1+\int_1^\infty t^{q-7/4}[\log(2+t)]^{3/4}\,dt\right)<\infty.
\]
Equation~\eqref{eq:typicallower} with $\eta=1/2$ gives the matching lower bound for large $n$; adjusting constants covers the remaining times.
\end{proof}

\begin{lemma}[Mean return in dimension seven]\label{lem:meanseven}
In $d=7$, $\E p_{2n}^H(o,o)\asymp n^{-4/5}$.
\end{lemma}
\begin{proof}
For $d=7$, the tail in \eqref{eq:unified-return-tail} is
$C[\log(2+t)]^{9/4}t^{-5/4}$, which is integrable. Thus
$\E p_{2n}\le Cn^{-4/5}$; Lemma~\ref{lem:returnlower} gives the reverse inequality.
\end{proof}

\subsection{Return fluctuations in dimensions five to seven}\label{sec:heat6}
For $d=5,6,7$, we tilt a middle block of the first Wilson walk. The unchanged prefix preserves local volume, while drift lowers the energy where the cutoff decreases. Write $g_G$ for the lattice Green function of Lemma~\ref{lem:lattice-walk}.

Let $\mathbb P_0$ be the ordinary joint law of Wilson's construction, with first walk $S$ started at $o$; its forest marginal is $\PP$. For an integer $N\ge1$ and $L\ge1$, obtain $\mathbb P_L$ from $\mathbb P_0$ by giving the increments of $S$ at $N\le k<2N$ normalized density $\exp(Le_1\cdot(S_{k+1}-S_k)/\sqrt N)$; all other increments and walks are unchanged. Write $\mathbb E_0,\mathbb E_L$ for the corresponding expectations, suppressing the dependence on $N$. For the last visit times $t_i=\max\{k\ge0:S_k=r_i\}$ of $S$ to $r_i\in\gamma=\LE(S)$, harmonically extend $f_i=F_N(t_i)$, where $F_N(t)=\min\{1,(2-t/N)_+\}$, to obtain $u$ on $H$.

\begin{lemma}[Tilted cutoff]\label{lem:tilted-cutoff}
Put $\kappa_5=\kappa_7=1$ and $\kappa_6=2$. For fixed $L\ge1$,
\begin{equation}\label{tiltl2}
 \mathbb E_0\left[(d\mathbb P_L/d\mathbb P_0)^2\right]\le e^{CL^2},
 \qquad
 \limsup_{N\to\infty}N^{(d-6)/2}\mathbb E_L\mathcal E_H(u)
       \le C L^{-\kappa_d}.
\end{equation}
There are constants $b,c,C,p>0$, independent of $L$, such that, for each fixed $L$ and all sufficiently large $N$, with probability at least $p$ under $\mathbb P_L$ the function $v=(u-1/4)_+$ has finite support, satisfies $\|v\|_\pi^2\le CN^2$, and is at least $1/4$ at $cN^2$ vertices of $H\cap Q_{b\sqrt N}$.
\end{lemma}
\begin{proof}
Throughout the proof, $L$ is fixed before $N\to\infty$; all constants are independent of $L$.

\emph{The energy under the tilted law.} The likelihood second moment is a product of one-step moment generating functions, bounded by $e^{CL^2}$. Under both laws the later Wilson walks depend on $S$ only through $\gamma$. Thus the conditional forest law given $(H,S)$ is the ordinary one given $(H,\gamma)$, for almost every $(H,S)$ under the tilt. Conditioning $N_{ij}$ first on $(H,S)$ and then on $S$ in Lemma~\ref{lem:ray-trace} gives
$\mathbb E_L[a(i,j)\mid S]=b_\gamma(i,j)\le C(1+|r_i-r_j|)^{2-d}$.
Since $r_i=S_{t_i}$ and the $t_i$ are strictly increasing, summing over all pairs of walk times in \eqref{eq:trace-form} gives
\begin{equation}\label{tiltenergy}
 \mathbb E_L\mathcal E_H(u)\le C\sum_{0\le s<t}
 (F_N(s)-F_N(t))^2\mathbb E_L(1+|S_t-S_s|)^{2-d}.
\end{equation}
For $N\ge L^2$, the maximum atom of $S_t-S_s$ is at most $C(t-s)^{-d/2}$ uniformly in $L$. At least half its increments have one homogeneous law, with uniformly positive step probabilities; Fourier inversion near the two periodicity points gives the bound, preserved by convolution with the remaining increments. Splitting at distance $\sqrt{t-s}$ bounds the Green weight by $C(t-s)^{(2-d)/2}$.

Fix $L$ and let $s/N\to a$, $t/N\to b>a$. Put $h=b-a$ and let $\delta$ be the length of $[a,b]\cap[1,2]$. The scaled increment converges to $Z_h+(L/d)\delta e_1$, where $Z_h$ is a centered Gaussian vector with covariance $(h/d)I_d$ and $I_d$ is the $d\times d$ identity matrix; the profile difference tends to $\delta$. For $\zeta>0$, the contribution from $|S_t-S_s|\le\zeta\sqrt N$ to the scaled Green weight is at most $C(N/(t-s))^{d/2}(\zeta^2+N^{-1})$. The rescaled integrand is dominated by the integrable function $C\min\{1,h^2\}h^{(2-d)/2}$ for $0\le a\le2$. The sums thus converge to integrals, bounded using the Gaussian estimate
\[
 \mathbb E[|Z_h+z|^{2-d}]\le C(h+|z|^2)^{(2-d)/2}
 \qquad(h>0,\ z\in\mathbb R^d).
\]
After removing the factor $N^{(6-d)/2}$, it remains to bound
\[
 \int_0^2\int_a^\infty
 \delta^2(h+L^2\delta^2)^{(2-d)/2}\,db\,da.
\]
If both endpoints lie in $[1,2]$, then $\delta=h$. The contribution is bounded by
$C\int_0^1 h^{(6-d)/2}(1+L^2h)^{(2-d)/2}\,dh$.
Splitting at $h=L^{-2}$ gives respectively $CL^{-3}\log(2L)$, $CL^{-2}$ and $CL^{-1}$ for $d=5,6,7$.
If exactly one endpoint lies outside, integrate first in $h-\delta$. The resulting bound is
\[
 C\int_0^1\delta^2(\delta+L^2\delta^2)^{(4-d)/2}\,d\delta
 \le
 \begin{cases}
 CL^{-1},&d=5,\\
 CL^{-2},&d=6,\\
 CL^{-3}\log(2L),&d=7.
 \end{cases}
\]
For endpoints on opposite sides of $[1,2]$, $\delta=1$ and integration in $b\ge2$ gives $CL^{4-d}$. Together these estimates prove \eqref{tiltl2}.

\emph{Mass near the origin.} Put $g_r(w)=\sum_{x\in Q_r}g_G(x,w)$ and $J_r=\sum_{i\ge0}g_r(S_i)$. Write $\mathbb E_{L,k,z}$ for expectation of the tilted walk started at $z\in\mathbb Z^d$ at absolute time $k\ge0$, with the same time-dependent increment laws. The atom bound remains valid after adding independent ordinary increments, uniformly over starting time and position. Hence
\[
 \sup_{k,z}\mathbb E_{L,k,z}\sum_{i\ge k}g_r(S_i)
 \le C\sum_{m\ge0}(m+1)\min\{1,r^d(m+1)^{-d/2}\}\le Cr^4.
\]
The Markov property gives $\mathbb E_LJ_r^2\le Cr^8$. The contribution to $\mathbb E_LJ_r$ from times $i\ge N/2$ is at most $Cr^dN^{(4-d)/2}$. If $q(w)$ is the probability of an independent lattice walk hitting $\gamma$, Wilson sampling gives
\[
 \mathbb P_L(x,y\in H\mid S)
 \le q(x)q(y)+C\sum_wg_G(x,w)g_G(y,w)q(w),\qquad
 q(w)\le\sum_ig_G(w,S_i).
\]
The second term accounts for the second walk meeting the first before it hits $\gamma$. Sum over $x,y\in Q_r$: the first term is controlled by $\mathbb E_LJ_r^2$, while $g_r(w)\le Cr^2$, $g_r(w)\le Cr^d|w|^{2-d}$ outside $Q_{2r}$, and $\sum_{w\in Q_h}\mathbb E_Lq(w)\le Ch^4$ control the second by dyadic summation. Thus $\mathbb E_L|H\cap Q_r|^2\le Cr^8$ uniformly in $L$.

Couple $S$ to an ordinary walk $S^0$ through time $N$. With probability tending to one, the ordinary walk has a cut time in $[N/2,3N/4]$. Indeed, on a two-sided ordinary walk the cut-time indicators are stationary and ergodic with positive mean, since independent half-walks have positive probability of disjoint ranges after the first step. The tilted continuation intersects the shared prefix up to $3N/4$ with probability at most $C\sum_{i\le3N/4,j\ge N}(j-i)^{-d/2}=O(N^{2-d/2})$. Off this last event, the cut is common to both walks, and their loop erasures agree through it. Compare the one-point marginals $\mathbb P_L(x\in H\mid S)$ and $\mathbb P_0(x\in H\mid S^0)$ by using the same independent test walk from $x$ to hit the respective rays. On the common-cut event a hit on the ordinary ray lost in the comparison must visit $S^0[N/2,\infty)$. Summing these probabilities costs at most $Cr^dN^{(4-d)/2}$. On the exceptional event, the ordinary sum is bounded by $J_r(S^0)$, whose second moment is at most $Cr^8$, so Cauchy--Schwarz gives a cost $o(r^4)$. The ordinary two-point estimate therefore yields
\[
 \mathbb E_L|H\cap Q_r|\ge c_0r^4-Cr^dN^{(4-d)/2}-o(r^4).
\]
Take $r=\lfloor b\sqrt N\rfloor$ with small fixed $b$. Paley--Zygmund gives $c_1r^4$ vertices with probability $p_1>0$. The conditional bush identity also bounds $\mathbb E_L\sum_{x\in H\cap Q_r}(1-u(x))$ by $Cr^dN^{(4-d)/2}$, since contributing attachments meet $S$ after time $N$. Decreasing $b$ and applying Markov leaves $cN^2$ vertices with $u\ge1/2$ with probability $3p_1/4$, uniformly in $L$ for sufficiently large $N$.

\emph{Finite support.} Let $\mathcal B$ be the union of the first $2N+1$ bushes. Lemmas~\ref{lem:past-estimates} and~\ref{lem:bush-domination} give conditional size tails $Cm^{-1/2}$. Truncation at $AN^2$, a union bound and Markov give $\mathbb P_L(|\mathcal B|>AN^2\mid S)\le C/\sqrt A$, without independence. For large fixed $A$, except on probability $p_1/8$, this event has conditional probability at most $1/8$ given $(H,S)$. Since $u(x)\le\mathbb P_L(x\in\mathcal B\mid H,S)$, every site with $u>1/4$ then contributes more than $1/8$ to
\[
 \mathbb E_L[|\mathcal B|\ind_{\{|\mathcal B|\le AN^2\}}\mid H,S]\le AN^2.
\]
Thus $|\{u>1/4\}|\le8AN^2$ and $\|v\|_\pi^2\le CN^2$. Intersecting with the mass event proves the claim with probability bounded below independently of $L$.
\end{proof}

The tilted cutoff gives large returns at many vertices; mass transport brings them to the root, and indistinguishability gives almost-sure fluctuations.
\begin{proposition}[Return fluctuations]\label{prop:return-fluctuations}
Almost surely,
\begin{gather}
 \limsup_n n^{4/3}p_{2n}^H(o,o)=\infty \qquad (d=5),\label{d5heatfluct}\\
 \limsup_n n p_{2n}^H(o,o)=\infty \qquad (d=6),\label{infinite-limsup}\\
 \limsup_n n^{4/5}p_{2n}^H(o,o)=\infty \qquad (d=7).\label{d7heatfluct}
\end{gather}
For $M\ge1$, with $\alpha=4/(d-2)$, the corresponding tail bounds are
\begin{equation}\label{tail6heat}
 \liminf_n\PP(n^\alpha p_{2n}^H(o,o)>M)\ge ce^{-CM^{\nu_d}},
 \qquad (\nu_5,\nu_6,\nu_7)=(3/2,1,5/2).
\end{equation}
\end{proposition}
\begin{proof}
\emph{Quantitative tails at the root.} Fix $D\ge1$, put $\alpha=4/(d-2)$ and take $t=\lfloor DN^{(d-2)/2}\rfloor$. On the mass and support event of Lemma~\ref{lem:tilted-cutoff}, the spectral inequality $|1-a^t|^2\le2t(1-a)$ for $a\in[-1,1]$ and integers $t\ge2$ gives $\|P_H^tv-v\|_\pi^2\le2t\mathcal E_H(v)$. Since $t\mathbb E_L\mathcal E_H(v)\le CDL^{-\kappa_d}N^2$ for sufficiently large $N$, choose $L=C_1D^{1/\kappa_d}$ with $C_1$ sufficiently large. Markov's inequality then makes the squared norm a sufficiently small multiple of $N^2$, jointly with the mass and support event, with probability bounded below. At least $cN^2$ vertices of $H\cap Q_{b\sqrt N}$ have $P_H^tv\ge1/8$. Reversibility gives
\[
 p_{2t}^H(x,x)\ge
 \frac{\pi(x)(P_H^tv(x))^2}{\|v\|_\pi^2}\ge cN^{-2}.
\]
Since $t^\alpha N^{-2}\asymp D^\alpha$, take $D=C_2M^{1/\alpha}$ with $C_2$ sufficiently large.

Write $E_t(x)=\{t^\alpha p_{2t}^{H_x}(x,x)>M\}$, and let $A$ be the event that at least $cN^2$ such vertices lie in $H\cap Q_r$, where $r=\lfloor b\sqrt N\rfloor$. For $Z=d\mathbb P_L/d\mathbb P_0$, Cauchy--Schwarz and \eqref{tiltl2} give $\mathbb P_0(A)\ge\mathbb P_L(A)^2/\mathbb E_0Z^2\ge ce^{-CL^2}$. The mass-transport calculation of Lemma~\ref{lem:returnlower} now yields
\[
 cN^2e^{-CL^2}
 \le\mathbb E_0\sum_{x\in H\cap Q_r}\ind_{E_t(x)}
 =\mathbb E_0[|H\cap Q_r|\ind_{E_t(o)}]
 \le Cr^4\mathbb P_0(E_t(o))^{1/2}.
\]
Since $r^4\asymp N^2$, the root probability is at least $ce^{-CL^2}$, with $L^2\asymp M^{2/(\alpha\kappa_d)}=M^{\nu_d}$. For any prescribed $n$, round $N$ upward so that $t\ge n$ and $t/n\to1$; monotonicity, with a slightly larger threshold, gives the full-time tails.

\emph{The almost-sure conclusion.} To transfer the lower tail to every root in a component, put $f_x=\ind_x/\pi(x)$, $h_{n,x}=P_H^nf_x$, and $q_n(x)=\|h_{n,x}\|_\pi^2=p_{2n}^H(x,x)/\pi(x)$. Then
\[
 \E\sum_{y\sim_H o}\|h_{n,o}-h_{n+1,y}\|_\pi^2
   =\E[p_{2n}-p_{2n+2}].
\]
Indeed, the two square terms have expectations $\E p_{2n}$ and $\E p_{2n+2}$ by mass transport, and $\sum_{y\sim_H o}f_y=\pi(o)P_H f_o$ identifies the cross term. Since $(\sqrt{q_n(y)}-\sqrt{q_{n+1}(y)})^2\le q_n(y)-q_{n+1}(y)$, the triangle inequality and another mass transport give
\[
 \E\sum_{y\sim_H o}(\sqrt{q_n(o)}-\sqrt{q_n(y)})^2
 \le4\E[p_{2n}-p_{2n+2}]
 \le Cn^{-1}\E p_{2\lfloor n/2\rfloor}.
\]
The last bound follows from $s^{2n}(1-s^2)\le Cn^{-1}s^{2\lfloor n/2\rfloor}$ on $[-1,1]$. Integrating \eqref{eq:unified-return-tail} gives $\E p_{2m}\le\E\sqrt{p_{2m}}\le Cm^{-\alpha/2}$. After multiplication by $n^\alpha$, the right-hand side is summable along $n=2^j$, since $\alpha<2$. Almost surely, the differences $n^{\alpha/2}(\sqrt{q_n(x)}-\sqrt{q_n(y)})$ therefore tend to zero along these times on every induced edge, by translation invariance and countability. In particular, $K(x)=\limsup_{j\to\infty}2^{j\alpha}q_{2^j}(x)$ is constant on each component.

For each integer $M\ge1$, the event $2^{j\alpha}p_{2^{j+1}}>4dM$ has probability bounded away from zero for large $j$ by \eqref{tail6heat}. Its probability of occurring infinitely often is at least the limsup of its probabilities. Since $\pi(o)\le2d$, this gives $\PP(K(o)>M)>0$. For a forest component $T'$, the condition that $K(x)>M$ for every vertex $x$ of $T'$ is an automorphism-invariant Borel property of the embedded tree: all the return probabilities are computed on $G[V(T')]$. Indistinguishability \cite[Theorem~1.1]{HN} and translation ergodicity \cite[Corollary~8.2]{BLPS} therefore give this property almost surely for every component. Intersecting over $M$ proves the assertion.
\end{proof}

The fixed-confidence scale and the preceding fluctuations rule out a deterministic normalization.
\begin{corollary}[Failure of deterministic normalization]\label{cor:return-fluctuation-consequences}
For $d=5,6,7$, no deterministic positive sequence $a_n$ admits eventual two-sided random-constant bounds for $p_{2n}^H(o,o)$.
\end{corollary}
\begin{proof}
Suppose such a sequence exists and write $\alpha=4/(d-2)$. Choose deterministic $K,N$ so that $K^{-1}a_n\le p_{2n}\le Ka_n$ for every $n\ge N$ with probability greater than $3/4$. Intersect at each $n$ with the fixed-confidence interval $[c_0n^{-\alpha},C_0n^{-\alpha}]$, also of probability greater than $3/4$. Then $c_0/K\le n^\alpha a_n\le KC_0$, and the assumed upper bound contradicts Proposition~\ref{prop:return-fluctuations}.
\end{proof}

\subsection{Six-dimensional return estimates}\label{sec:d6}
Throughout this subsection $d=6$. We sharpen confinement for the pointwise lower bound, use the Green function and the box resistance for the upper and cumulative bounds, and integrate the cut tail for fractional moments.
\begin{proposition}[Six-dimensional return estimates]\label{prop:d6-return-bounds}
For every $\varepsilon>0$, almost surely eventually,
\begin{equation}\label{explicitheat6}
 \frac{c_\varepsilon}{n(\log n)^{5+\varepsilon}}
 \le p_{2n}^H(o,o)\le
 \frac{C_\varepsilon\log n(\log\log n)^{1+\varepsilon}}n.
\end{equation}
Moreover, for $n\ge1$,
\begin{align}
 c/n\le\mathbb E p_{2n}^H(o,o)&\le C\lp_n/n,\label{d6mean}\\
 \mathbb E\sum_{k=0}^np_{2k}^H(o,o)&\asymp\lp_n,\label{integratedheat}\\
 c_qn^{-q}\le\mathbb E[p_{2n}^H(o,o)^q]&\le C_qn^{-q}\qquad(0<q<1).\label{d6fractional}
\end{align}
\end{proposition}
\begin{proof}
\emph{Pointwise lower bound.}
Fix a small failure probability $\delta$. In \eqref{eq:confinement-general}, take $N=\lfloor aR^2\rfloor$, $r=\lfloor bR\rfloor$ and fix $p>2$ with conjugate exponent $p'=p/(p-1)$. Choose
\[
 v=c[\log(2/\delta)]^{-3},\qquad
 a=c_p(v\delta)^{p'},\qquad b^2=cav\delta.
\]
The expected exterior mass is $O(\delta)$ whenever $n\le c\delta vr^4$. Markov's inequality makes the quenched exterior mass at most $1/2$ outside probability $C\delta$. The least admissible $R$ satisfies
$R^4\le C_pn(v\delta)^{-(2p'+3)}$.
The volume bound $V_{2R}\le C_p\delta^{-1/p}R^4$, which fails with probability at most $\delta$, and the confinement-to-return inequality give
\[
 p_{2n}\ge c_pn^{-1}\delta^{1/p}(v\delta)^{2p'+3}
\]
outside probability $C\delta$. Set $n=2^j$ and $\delta=(1+j)^{-1-\eta}$. Since $1/p+2p'+3\to5$, choose $p$ large and then $\eta$ small. The factor $v$ contributes only a power of $\log\log n$, which is absorbed in the remaining $\varepsilon$ slack. Borel--Cantelli and monotonicity prove the lower half of \eqref{explicitheat6}.

\emph{Green-function bounds.}
In $2n$ steps the walk cannot leave $Q_{2n}$. Monotonicity of even returns and the killed Green identity $G_{H\cap Q_R}(o,o)=\pi(o)\mathcal R_R$ from the proof of Lemma~\ref{lem:terminal} therefore give
\begin{equation}\label{finiteprop}
 (n+1)p_{2n}\le\sum_{k=0}^np_{2k}\le\pi(o)\mathcal R_{2n}.
\end{equation}
By \eqref{eq:six-mean-resistance}, $\E\mathcal R_{2n}\le C\lp_n$, which gives the upper halves of \eqref{d6mean} and \eqref{integratedheat}. The lower halves follow from \eqref{eq:typicallower}, which gives $\E p_{2k}\ge c/k$. For the almost-sure upper bound, put $R_k=\lfloor e^{2^k}\rfloor$. Markov's inequality and \eqref{eq:six-mean-resistance} give $\PP(\mathcal R_{R_k}>k^{1+\varepsilon}\log R_k)\le Ck^{-1-\varepsilon}$, which is summable. By Borel--Cantelli and monotonicity of $\mathcal R_R$ in $R$, eventually $\mathcal R_R\le\mathcal R_{R_{k+1}}\le C(k+1)^{1+\varepsilon}\log R_{k+1}$ for $R_k\le R<R_{k+1}$. Since $\log R_{k+1}\le3\log R$ and $k\le C\log\log R$, this gives $\mathcal R_R\le C_\varepsilon\log R(\log\log R)^{1+\varepsilon}$ eventually. With \eqref{finiteprop}, this proves the upper half of \eqref{explicitheat6}.

\emph{Fractional moments.}
The unified cut tail \eqref{eq:unified-return-tail}, specialized to $d=6$, is
\begin{equation}\label{d6nptailupper}
 \mathbb P(np_{2n}>Ct)\le C[\log(2+t)]^{3/2}/t
 \qquad(n,t\ge1).
\end{equation}
At dyadic times $n=2^j$, summability of this tail bound at $t=j(\log j)^{5/2+\varepsilon}$ would give an almost-sure upper bound with power $5/2+\varepsilon$ on $\log\log n$. The resistance estimate above gives power $1+\varepsilon$. Integrating against $qt^{q-1}$ proves the upper bound in \eqref{d6fractional} for $0<q<1$; the lower bound follows from \eqref{eq:typicallower}.
\end{proof}

\section{Repeated contacts in dimension eight}\label{sec:contact-sharp}
We prove the logarithmic improvements in dimension eight by using reciprocal contact weights to count each neighboring component once. Local Wilson comparisons yield joint contact estimates; a second moment bounds the reciprocal count. We then transfer this bound to the origin-wired law and apply the deletion identities to resistance and returns.

Throughout this section, $d=8$ and $e$ is a fixed unit lattice vector. We use the bracket convention $\langle x\rangle=1+|x|$ from Section~\ref{sec:notation} and the lattice Green function $g_G$ from Lemma~\ref{lem:lattice-walk}. Recall that $T_x$ is the ordinary WUSF component of $x$, and put
\[
 p(x)=(T_x,T_{x+e}),\qquad \mathcal D_x=\{T_x\ne T_{x+e}\},\qquad
 N_r(x)=\sum_{w\in Q_r(x)}\ind_{\{p(w)=p(x)\}}.
\]
The order of the two components is retained in $p(x)$. Radii below are rounded to integers.

\begin{proposition}[Reciprocal contact count]\label{prop:contact-multiplicity}
For $r\ge1$,
\begin{equation}\label{eq:sharp-contact}
 q(r):=\E[\ind_{\mathcal D_0}/N_r(0)]\asymp\frac1{\log(2+r)}.
\end{equation}
\end{proposition}
The first moment gives the lower bound; two- and three-contact comparisons control the decrease of $q$ between scales.

\subsection{Local comparison of Wilson forests}\label{sec:contact-gain}
We first control occupation and local loop erasures for walks conditioned on prescribed visits.
\begin{lemma}[Occupation of a conditioned walk]\label{lem:bridge-occupation}
Let $S$ be lattice simple random walk under $P_v^G$ and $\tau_h=\inf\{n\ge0:S_n=h\}$. Conditioning on $\tau_h<\infty$ leaves ordinary continuation after that hit, and
\begin{equation}\label{eq:bridge-occupation}
\begin{aligned}
 E_v^G\!\left[\left.\sum_{n\ge0}\ind_{\{S_n=w\}}\right|\tau_h<\infty\right]
 &\le\frac{g_G(v,w)g_G(w,h)}{g_G(v,h)}+g_G(h,w)\\
 &\le C\{g_G(v,w)+g_G(h,w)\}.
\end{aligned}
\end{equation}
Stopping at the prescribed hit only decreases the occupation.
For two prescribed vertices visited successively, sum this bound over the stages; coinciding vertices give a zero-length stage.
\end{lemma}
\begin{proof}
The Doob transform \cite[Proposition~11.1.1]{LL} gives the first term before the hit, and continuation contributes $g_G(h,w)$. Apply \eqref{eq:lattice-green} and the triangle inequality, then the Markov property for successive hits.
\end{proof}
The next coupling controls the effect of distant roots.

\begin{lemma}[Locality under additional wiring]\label{lem:contact-locality}
For any absorbing set $A$, Wilson forests rooted at $A\cup\{\infty\}$ and additionally at $Q_T^c$ can be coupled to agree on the internal edges of $Q_s$ with failure probability at most $C(s/T)^6$, provided $4s<T$. An ordinary walk has a cut time in $[N,2N]$ with failure probability at most $CN^{-1/2}$.
\end{lemma}
\begin{proof}
Let $\mathcal G_A$ be the inverse of the Dirichlet graph Laplacian $Lf(v)=\sum_{u\sim v}(f(v)-f(u))$, with functions set to zero on $A$. The kernel $K=\mathcal G_A-\mathcal G_{A\cup Q_T^c}$ is positive semidefinite, and the strong Markov property at exit gives $K(v,v)\le CT^{-6}$ on $Q_{T/2}$. Write $K(u,v)=(f(u),f(v))$ in a Hilbert space; $f$ vanishes on $A$ and is harmonic on $Q_{T-1}\setminus A$. Summation by parts against $\eta^2f$, with $\eta=1$ on $Q_{s+1}$, supported in $Q_{2s+2}$ and having edge differences at most $C/s$, gives
\[
 \sum_{\{u,v\}\subset Q_s,\,u\sim v}\|f(u)-f(v)\|^2
 \le Cs^{-2}\sum_{Q_{2s+2}}\|f(u)\|^2\le C(s/T)^6.
\]
Use $2ab\le a^2/2+2b^2$, absorbing bounded marginal cases into $C$. For $e=\{u,v\}$ the transfer-current theorem \cite[Ch.~4, (4.5)]{LP}, applied after contracting the roots and then taking the wired limit, identifies the decrease in $\PP(e\in F)$ with $K(u,u)+K(v,v)-2K(u,v)$. This is precisely $\|f(u)-f(v)\|^2$. Under the inclusion coupling of Lemma~\ref{lem:wilson-facts}, the sum therefore bounds the probability of any discrepant internal edge. Wiring the exteriors of disjoint boxes makes the remaining forests independent; the same bound applies conditionally whenever previously constructed paths agree in the larger boxes.

For cut times, choose $c\sqrt N$ candidates with disjoint $\lfloor\sqrt N\rfloor$-step windows. Each event that the past and future inside one window are disjoint depends only on that window's increments. These events are independent and each has probability bounded below by the probability that two independent half-infinite walks do not meet after their initial separation, which is positive in $d=8$. An intersection crossing a candidate from outside its window has probability at most $\sum_{j\ge\lfloor\sqrt N\rfloor}(j+1)p_j^G(0,0)\le CN^{-1}$. A union bound completes the proof.
\end{proof}

This allows a local comparison with ordinary WUSF under distant conditioning.
\begin{lemma}[Distant targets]\label{lem:contact-distant}
Run Wilson rooted at $b$ and infinity, with $|b|\ge D$. Start from $0$ with a walk conditioned to hit $b$, then from $e$ with an ordinary walk or one conditioned to hit $h$, $|h|\ge D$. An optional third walk from $z$, $|z|\ge D$, may also be conditioned to hit $h$. The driving walks are independent and continue ordinarily after their prescribed hits. Subsequent walks are ordinary. The internal edges of $Q_{\lfloor D^{1/1000}\rfloor}$ have total variation distance at most $CD^{-c}$ from ordinary WUSF, uniformly also for coinciding targets and $z=h$.
\end{lemma}
\begin{proof}[Proof sketch]
For a walk started near the center, the Green function ratio couples its path up to exit from $Q_P$ to an ordinary walk with error $O(P/D)$, by \eqref{eq:green-expansion}. If its start and target are both at distance at least $D$, the probability of entering $Q_P$ is $O((P/D)^6)$ by the same ratio and the ordinary entrance estimate. After exit from $Q_P$, the probability of returning to $Q_u$ is $O((u/P)^6+(u/D)^6)$ for $4u\le P\le D/4$.

To preserve the erasure in $Q_U$, couple the walks up to exit from $Q_{U\ell^{4M}}$, where $U\ell^{4M}\le D/4$ and $\ell^M\ge4$. Put $N=\lceil U^2\ell^{2M}\rceil$. Lemma~\ref{lem:contact-locality} gives a cut time in $[N,2N]$. The heat kernel bound keeps the prefix in $Q_{U\ell^{2M}}$ and, by last exit, bounds the chance of a visit to $Q_U$ after $N$ by $CU^6N^{-3}$. The return estimate above now gives a common cut with total error
\[
 \varepsilon_D(U)=C\big[(U\ell^M)^{-1}+\ell^{-6M}+U\ell^{4M}/D\big].
\]
A later stop preserves the local erasure. For two successive walks, protect the first branch in $Q_{U_1}$ with $U_1\ge U\ell^{2M}$, so any hit before the second cut occurs on the same branch in both constructions. Wiring $Q_U^c$ then couples their completions on $Q_s$ with additional error $C(s/U)^6$. Taking $s=\lfloor D^{1/1000}\rfloor$, $U=\ell=D^{1/100}$, $M=1$ and $U_1=D^{1/25}$ gives the claim.
\end{proof}

For the lower comparison we also need to keep the two components distinct.
\begin{lemma}[Two prescribed component connections {\cite[Lemma~4.21]{HP}}]\label{lem:prescribed-contact-pairs}
Let $x_1,x_2,y_1,y_2$ be distinct vertices. If four independent walks from these vertices are mutually disjoint with probability at least $\epsilon>0$, then
\[
 \PP(T_{x_1}=T_{y_1},\ T_{x_2}=T_{y_2},\ T_{x_1}\ne T_{x_2})
 \ge c_\epsilon(1+\diam_{\rm ext}\{x_1,x_2,y_1,y_2\})^{-8}.
\]
\end{lemma}
We verify the nonintersection hypothesis when applying this lemma.

\subsection{Probabilities of separated contacts}\label{sec:contact-localization}
We localize the contact counts before comparing distant boxes. Let $N_{r,S}(x)$ count $w\in Q_r(x)$ for which $w$ connects to $x$ and $w+e$ to $x+e$ in the restricted forest $F[Q_S(x)]$. Set $\phi_{r,S}(x)=0$ if $x$ and $x+e$ connect in this restricted forest, and $\phi_{r,S}(x)=1/N_{r,S}(x)$ otherwise; $N_{r,S}(x)\ge1$ for $S\ge r+2$.

\begin{lemma}[Local contact counts]\label{lem:contact-local-count}
For $S\ge4r+10$,
\begin{equation}\label{eq:near-weight}
 q(r)\le\E[\ind_{\mathcal D_0}\phi_{r,S}(0)]\le\E\phi_{r,S}(0)
 \le q(r)+C(r/S)^4+CS^{-4}.
\end{equation}
Moreover, $q(r)\ge c/\log(2+r)$.
\end{lemma}
\begin{proof}
Two vertices of $Q_{r+1}$ connect outside $Q_S$ with probability at most $CS^{-4}$: a Wilson intersection outside $Q_{S/2}$ costs $\sum_{u\notin Q_{S/2}}g_G(v,u)g_G(v',u)\le CS^{-4}$, and one inside after an earlier exit costs $CS^{-6}\cdot CS^2$. An omitted contact requires this in one component and a connection of probability at most $C\langle w\rangle^{-4}$ in the other. Run Wilson by component and discard intergroup avoidance to multiply the independent walk bounds. Summing over $w\in Q_r$ gives $C(r/S)^4$; global coalescence without local coalescence costs $CS^{-4}$, proving \eqref{eq:near-weight}. The same independent walk bounds give $\PP(\mathcal D_0,\,p(w)=p(0))\le C\langle w\rangle^{-8}$. Summing over $w$ gives $\E[\ind_{\mathcal D_0}N_r(0)]\le C\log(2+r)$. Since two walks from $0,e$ have positive probability of disjoint ranges, $\PP(\mathcal D_0)>0$, and Cauchy--Schwarz proves the lower bound.
\end{proof}

For the second moment we compare two and three contacts. Let $I$ mean that $p(x_i)$ agree and have distinct entries, for $m=2,3$ distinct centers. Put $t=\max_{i,j}\langle x_i-x_j\rangle$, $\Delta=\min_{i\ne j}\langle x_i-x_j\rangle$, $\ell=\log(2+t)$, and $\mathcal T=t$ for $m=2$, $\mathcal T=t\Delta$ for $m=3$. Assume $\Delta\ge t/\ell^{C_*}$ for a fixed $C_*$. Write $k(x,y)=\langle x-y\rangle^{-4}$, and for lattice kernels use $(f_1*f_2)(x,y)=\sum_z f_1(x,z)f_2(z,y)$.

\begin{lemma}[Separated contact comparison]\label{lem:contact-tests}
For every $B>1$ there is $L_B$ such that, if $s\le\Delta/\ell^{L_B}$ and $\varphi_i\in[0,1]$ depend on internal edges of $Q_s(x_i)$, then, for any $J\subseteq\{1,\ldots,m\}$,
\begin{equation}\label{eq:decor-up}
 \E[\ind_I\prod_{i\in J}\varphi_i]
 \le C\mathcal T^{-8}\{\prod_{i\in J}\E\varphi_i+\ell^{-B}\}.
\end{equation}
For two centers,
\begin{equation}\label{eq:decor-down}
 \E[\ind_I\varphi_1\varphi_2]
 \ge ct^{-8}\{\E[\ind_{\mathcal D_{x_1}}\varphi_1]
                  \E[\ind_{\mathcal D_{x_2}}\varphi_2]-C\ell^{-B}\}.
\end{equation}
In the two-center case, for either center $c$ and $4r+10\le S\le t/\ell^{2B+10}$,
\begin{equation}\label{eq:weighted-near}
 \PP(I,\ N_{r,S}(c)<N_r(c))
 \le Ct^{-8}\{(r/S)^4+\ell^{-B}\}.
\end{equation}
The exponent $L_B$ depends only on $B$. The constants and the threshold for $t$ may depend on $C_*$.
\end{lemma}
\begin{proof}[Proof sketch]
On $I$ the walks from the $x_i$ form a connected intersection graph, as do those from the $x_i+e$. Mark one meeting in each group for two centers, or two ordered meetings for three. Dropping avoidance between the groups gives independent walk diagrams. For three walks the weight is a sum of terms
$g_G(a,h)g_G(b,h)g_G(h,h')g_G(c,h')$.
Convolution bounds the weight of each group by $C\mathcal T^{-4}$ \cite[Lemma~2.12]{BKPS}. A mark within $\delta\Delta$ of a center gains a factor $C\delta^2$ by summing the Green kernel in that box. We may therefore discard marks within $D=\Delta/\ell^{B+2}$ at cost $C\mathcal T^{-8}\ell^{-B}$.

For each remaining diagram, the independent walks are conditioned to visit their marks in the prescribed order and continue ordinarily after the last hit. Apply the entrance bound at each stage. At each center the common-cut construction in Lemma~\ref{lem:contact-distant} preserves the two local branches. For example, take $M=B+10$, $L_B=60M$, $U=\Delta/\ell^{30M}$ and $U_1=U\ell^{6M}$. The cut, entrance and wiring errors together are $O(\ell^{-B})$. Wiring outside the disjoint boxes makes the completions independent, giving \eqref{eq:decor-up}.

For \eqref{eq:decor-down}, expose four prefixes up to their exits from boxes of radius $a=t/\ell^{20M}$. Continue them independently in the two boxes' comparison forests. From the same exits start fresh walks and let $E$ be the prescribed pairing across the boxes in distinct trees. When the two exits in each box are separated by $a/\ell^M$, Lemma~\ref{lem:prescribed-contact-pairs} gives $\PP(E\mid\text{prefixes})\ge ct^{-8}$, since four walks from these exits avoid one another with probability bounded below.

Use the fresh walks also to continue the original prefixes, running both left starts before both right starts. A common cut protects the local erasure in $Q_\rho$, with $\rho=t/\ell^{35M}$. A second cut occurs before exit from $Q_b$, with $b=t/\ell^{10M}$, and its tail avoids $Q_{2a}$, so restoring the prefix preserves that tail. On $E$, if the starts in each comparison forest belong to distinct trees and there are no unintended intersections of prefixes and continuations, the right branches attach outside both $b$-boxes at the exit forest's stopping vertices. To bound the exceptions jointly with $E$, mark the two cross-box meetings. The entrance and occupation bounds of Lemmas~\ref{lem:bridge-occupation} and~\ref{lem:contact-distant} give total cost $O(t^{-8}\ell^{-B})$ for poor exit separation, near marks, failed cuts and unintended intersections. Given the exposed walks, wiring outside the $\rho$-boxes couples the completions in $Q_s$ with error $O((s/\rho)^6)$ by Lemma~\ref{lem:contact-locality}. The comparison forests are independent before imposing $E$, whose conditional probability given the prefixes is at most $Ct^{-8}$. Integrating the conditional lower bound and subtracting these errors gives \eqref{eq:decor-down}.

Finally, start Wilson at the tested center $c$ and its translate. With the distant marks fixed, an omitted contact at $w$ requires a hit on the far branch or a local connection leaving $Q_S(c)$. By \eqref{eq:bridge-occupation} these cost $O(D^{-4}+S^{-4})$ together. The other group costs $O(k(c,w)+D^{-4})$. Summing their product over $Q_r(c)$ gives $O((r/S)^4)$. Restoring the diagram weight and the discarded marks gives \eqref{eq:weighted-near}.
\end{proof}

\subsection{The reciprocal count and the origin-wired law}
A new box with at least as many contacts doubles the count. We find such a box by a weighted second moment.
\begin{proof}[Proof of Proposition~\ref{prop:contact-multiplicity}]
The lower bound is Lemma~\ref{lem:contact-local-count}.

\emph{Doubling the contact count.} Fix $B=12$, $L>2L_B+2B+20$, and put
$\lambda=\log(2+r)$, $S=r\lambda^L$, $r_-=r/\lambda^{2L}$,
$S_-=r/\lambda^L$, $T=r\lambda^{6L}$ and
$W=\{z:r\lambda^{4L}\le|z|\le T/3\}$.
Set $I_z=\{p(z)=p(0)\}\cap \mathcal D_0$, $A=N_r(0)$, $B_z=N_r(z)$ and
$m_0=\log(T/r)\asymp\log\lambda$. Equations
\eqref{eq:decor-down}, \eqref{eq:near-weight} and \eqref{eq:weighted-near} imply
\begin{equation}\label{eq:rare-weight-lower}
 \E[\ind_{I_z}/(AB_z)]\ge c|z|^{-8}q(r)^2,
\end{equation}
using $q(r)\ge c/\lambda$ to absorb errors. With $W_0=\ind_{\mathcal D_0}/A$, define
\[
 X=\sum_{z\in W}\ind_{I_z}\ind_{\{B_z\ge A\}}/B_z.
\]
On $X>0$, two disjoint $r$-boxes in $Q_T$ contain at least $2A$ contacts, so $q(r)-q(T)\ge\frac12\E[W_0\ind_{\{X>0\}}]$. The point reflection $u\mapsto z+e-u$ takes the oriented edge at $w$ to the reversed edge at $z-w$. It preserves $I_z$ and exchanges $A$ and $B_z$, so the ordering retains at least half the expectation in \eqref{eq:rare-weight-lower}. We obtain
\begin{gather}
 \E[W_0X]\ge cm_0q(r)^2,\label{eq:rare-first-moment}\\
 \E[W_0X^2]\le C\{m_0q(r_-)^2+m_0^2q(r_-)^3\}.\label{eq:rare-second-moment}
\end{gather}
\emph{Second moment.} Discard the ordering. For near-diagonal pairs, first use the global counts: covering $Q_r(z)$ by boundedly many subcubes of diameter $r$ gives
\[
 \sum_{\substack{w\in Q_r(z)\\p(w)=p(z)}}N_r(w)^{-1}\le C_8.
\]
Then bound the two remaining reciprocals by $\phi_{r_-,S_-}$ at their respective centers, each of mean at most $Cq(r_-)$. The two-center bound \eqref{eq:decor-up} gives the first term of \eqref{eq:rare-second-moment}. For the other pairs, $\Delta\ge r$, $t\le CT$ and $\log(2+t)\asymp\lambda$; take $C_*=6L+2$. Since $L>2L_B$, localizing all three reciprocals in boxes of radius $S_-$ satisfies the hypotheses of \eqref{eq:decor-up}. The connection weight is at most $C[|z|\min(\langle w\rangle,\langle w-z\rangle)]^{-8}$; its sum over $z,w$ is $O(m_0^2)$, giving the second term. The errors are absorbed by $q(r_-)\ge c/\lambda$.

Weighted Cauchy--Schwarz gives
\[
 q(r)-q(T)\ge\frac{cm_0q(r)^4}{q(r_-)^2[1+m_0q(r_-)]}.
\]
\emph{Iteration over scales.} Write $q_k=q(\lfloor e^k\rfloor)$. For $K/2\le k\le K$ choose
$m=\lceil C_L\log K\rceil$ covering the logarithms of $r_-,T$.
Monotonicity gives $q_{k-m}-q_{k+m}\ge cmq_{k+m}^4/
\{q_{k-m}^2[1+mq_{k-m}]\}$. If $q_{k+m}\ge q_{k-m}/2$, division by $q_{k-m}q_{k+m}$ gives a reciprocal increment at least $cm/(1+mq_{k-m})$. Otherwise the increment is at least $1/q_{k-m}$. Since $q_{k-m}\le1$, both bounds are bounded below by a positive constant. Summing over disjoint intervals $[k-m,k+m]\subset[K/2,K]$ gives $q_K\le Cm/K$. Apply this first bound also at $\lfloor K/3\rfloor$. Monotonicity now gives $mq_{k-m}=O(m^2/K)=o(1)$ throughout the block. Both reciprocal increments are then at least $cm$. Summing again over order $K/m$ intervals gives $q_K\le C/K$.
\end{proof}

We transfer this estimate to the origin-wired law, where reciprocal weights count all contacts with each component once.
\begin{lemma}[Origin-wired contact bounds]\label{lem:wired-contact}
Under the conditioned Wilson laws of Lemma~\ref{lem:contact-distant}, the inverse full contact multiplicity of the central components, set to zero when they coalesce, has mean at most
\begin{equation}\label{eq:sharp-contact-factor}
 C/\log(2+D).
\end{equation}
In the origin-wired forest, let $K_z$ be the component of $z\notin A_o$. Then
\begin{gather}
 \PP^{o,\mathrm w}(z\notin A_o,\ K_z\sim A_o)
 \le C\langle z\rangle^{-2}/\log(2+|z|),\label{eq:d8-two-site}\\
 \E\sum_{z\in H\cap Q_R}\res^H(o,z)
 \le CR^6/\log(2+R).\label{eq:d8-pair-sum-explicit}
\end{gather}
\end{lemma}
\begin{proof}
\emph{Conditioned reciprocal weights.} Let $\E_\star$ denote expectation under any one of the conditioned Wilson laws in Lemma~\ref{lem:contact-distant}. For $s=\lfloor D^{1/1000}\rfloor$, $r=\lfloor\sqrt s\rfloor$ and large $D$, Lemma~\ref{lem:contact-distant} and \eqref{eq:near-weight}, \eqref{eq:sharp-contact} give
\begin{equation}\label{eq:local-contact-factor}
 \E_\star[\phi_{r,s}(0)]\le q(r)+C(r/s)^4+Cs^{-4}+CD^{-c}
 \le C/\log(2+D).
\end{equation}
The full reciprocal is at most $\phi_{r,s}(0)$, proving \eqref{eq:sharp-contact-factor}; absorb bounded $D$ into $C$.

\emph{Summing over contacts.} Conditional on $x\in A_o$, run Wilson from $x,y=x+e,z$. On $y\notin A_o$, $z\in K_y$, the independent driving walks from $y,z$ meet. Bound this event by summing over common vertices $h$. For each $h$, condition both walks to hit it (even if Wilson stops a walk earlier). The product probability is at most $Cg_G(y,h)g_G(z,h)$. Apply the comparison translated by $x$, with
$D=\min(|x|,|z-x|,|h-x|)$ whenever $D$ is large, and use the bound one on the reciprocal otherwise. With the integrand defined as zero for $y\in A_o$, this gives
\[
 \E^{o,\mathrm w}\!\left[\left.\frac{\ind_{\{y\notin A_o,\ z\in K_y\}}}{m(A_o,K_y)}
 \right|x\in A_o\right]
 \le\frac{C\langle z-x\rangle^{-4}}
 {\log(2+\min(|x|,|z-x|))}.
\]
For $q=\min(|x|,|z-x|)<4$, the unweighted convolution $g_G*g_G\le Ck$ suffices, including $x=o$ and coinciding starts. For $q\ge4$, meetings with $|h-x|\le\sqrt q$ cost at most $Cq\langle z-x\rangle^{-6}$. Outside this box the logarithm has order $\log(2+q)$, and we again use $g_G*g_G\le Ck$.

Multiply by $\PP^{o,\mathrm w}(x\in A_o)=g_G(x,o)/g_G(o,o)$ and sum over $x,e$; the reciprocal contact weights for $K_z$ sum to one. The unweighted convolution is $C\langle z\rangle^{-2}$. Writing $B=|z|$, the regions $|x|\le\sqrt B$ and $|z-x|\le\sqrt B$ contribute $CB^{-3}$ and $CB^{-4}$; elsewhere the logarithm has order $\log B$. This proves \eqref{eq:d8-two-site}, and two-point deletion gives \eqref{eq:d8-pair-sum-explicit}.
\end{proof}

\subsection{Resistance and return probabilities}\label{sec:d8obs}\label{sec:qual8}
Reciprocal contact weights now improve the resistance and return bounds by a logarithm.
\begin{proposition}\label{prop:d8-sharp}
At each fixed confidence level and in mean, the respective scales of $\mathcal R_R$ and $p_{2n}^H(o,o)$ are
\begin{equation}\label{conf8}
 R^2/\log(2+R),\qquad[n\log(2+n)]^{-2/3}.
\end{equation}
The normalized return probabilities also have uniformly bounded $q$th moments for $0<q<3/2$.
\end{proposition}
\begin{proof}
\emph{Lower bounds.} The return lower bound is \eqref{eq:typicallower}. For the resistance, apply Lemma~\ref{sec:spatial} with $D=H\cap Q_R$ and $N=\lfloor aR^2\rfloor$, where $0<a<1$. On $\{q_D\le1/4\}$, the function $(2u_N-1)_+$ equals one at $o$, vanishes outside $Q_R$ and has energy at most $4\mathcal E_H(u_N)$, so Dirichlet's principle gives $\mathcal R_R\ge1/(4\mathcal E_H(u_N))$. By \eqref{energies}, $\E\mathcal E_H(u_N)\le C\log N/N$. For $0<c'<1$, \eqref{supportfail} and Markov's inequality give, for large $R$,
\[
 \PP\big(\mathcal R_R<c'R^2/\log R\big)\le C(a+e^{-c/a})+\frac{4c'R^2}{\log R}\cdot\frac{C\log N}{N}\le C(a+e^{-c/a})+\frac{Cc'}a.
\]
Choosing $a$ small and then $c'$ small gives the fixed-confidence lower bound, and hence the mean lower bound.

\emph{Resistance.} Write $A=A_o$. By \eqref{terminal},
$\E\mathcal R_R=r_G\E^{o,\mathrm w}\sum_{K\sim A}P_o^{G[A\cup K]}(X_{\tau_R}\in K)$.
Each positive summand has a contact $(v,v+e)$ with $v\in A\cap Q_R$, $v+e\in K$. For this direction put $P_{K,e}=\{x\in A:x+e\in K\}$, $s=\lceil\sqrt R\rceil$, $N_s^e(x)=|P_{K,e}\cap Q_s(x)|$. Then
\begin{equation}\label{localpack8}
 \sum_{x\in P_{K,e}\cap Q_{2R}}1/N_s^e(x)\ge1/2.
\end{equation}
Use the boxes $Q_{js}(v)$, $0\le j\le J=\lceil\log_2 |Q_{2R}|\rceil+1$. For large $R$, $Js\le R$, and all these boxes lie in $Q_{2R}$. The first contains $v$, so the counts cannot double at every step. At a step where doubling fails, every contact of the smaller box has its $s$-box inside the larger box. Summing its reciprocal count over the smaller box gives at least $1/2$, proving \eqref{localpack8}.

For $x\in A$, give $(x,x+e)$ weight zero if $x+e\in A$, and weight $1/N_s^e(x)$ otherwise, taking $K$ to be the tree of $x+e$. For large $D=|x|$ and $u=\lfloor D^{1/1000}\rfloor$, we have $\lfloor\sqrt u\rfloor\le s$ on $Q_{2R}$. The weight is therefore bounded by $\phi_{\lfloor\sqrt u\rfloor,u}(x)$, of conditional mean at most $C/\log(2+D)$ given $x\in A$ by \eqref{eq:local-contact-factor}. Bounded $D$ contribute $O(1)$. Summing \eqref{localpack8} and using $\PP^{o,\mathrm w}(x\in A)=g_G(x,o)/g_G(o,o)$ gives
$\E\mathcal R_R\le C\sum_{x\in Q_{2R}}\langle x\rangle^{-6}/\log(2+|x|)\le CR^2/\log(2+R)$.

\emph{Return probabilities.} Apply \eqref{d8av} with $V_r=|H\cap Q_r|$ and
$Y_r=\sum_{z\in H\cap Q_r}\res^H(o,z)$. Lemma~\ref{lem:forest-volumes} and \eqref{eq:d8-pair-sum-explicit} give
$\PP(V_r<br^4)\le Ce^{-cb^{-1/3}}$ and $\E Y_r\le Cr^6/\log(2+r)$.
Put $R_0=(n\log(2+n))^{1/6}$. For $0<a<1/4$, $0<b_0<1-4a$ and $1\le\lambda\le R_0^4$, take
$r=\lfloor R_0\lambda^{-a}\rfloor$, $b=\lambda^{-b_0}$.
Then $\log(2+r)\asymp_a\log(2+n)$, and $4a+b_0<1$ makes the volume term at most half of $C_0\lambda/R_0^4$ on $V_r\ge br^4$. Markov's inequality on the other term gives
\[
 \PP(p_{2n}^H(o,o)>C_0\lambda/R_0^4)
 \le Ce^{-c\lambda^{b_0/3}}+
 \frac{CR_0^4r^2}{\lambda nb\log(2+r)}
 \le C_ae^{-c\lambda^{b_0/3}}+C_a\lambda^{-(1+2a-b_0)}.
\]
For larger $\lambda$ the event is empty. Taking $a=b_0=1/8$ gives exponent $9/8$, hence the mean and fixed-confidence upper bounds. Letting $a\uparrow1/4$, $b_0\downarrow0$ and integrating gives all moments below $3/2$.
\end{proof}

\section{Open questions}\label{sec:questions}
\begin{enumerate}[leftmargin=*]
\item In dimension six, is $\E p_{2n}^H(o,o)\asymp n^{-1}$? All strictly smaller positive moments have the predicted order, and the cumulative mean grows as $\log n$. Is $\sum_{k\le n}p_{2k}^H(o,o)\ge c\log n$ eventually almost surely, and how small can $np_{2n}^H(o,o)$ be along subsequences?
\item In dimension five, is $\E p_{2n}^H(o,o)\asymp n^{-4/3}$? Our moment estimate $\E[(p_{2n}^H(o,o))^q]\asymp n^{-4q/3}$ covers $0<q<3/4$.
\item In dimension eight, for which $q>0$ are the $q$th moments of $(n\log n)^{2/3}p_{2n}^H(o,o)$ bounded uniformly in $n\ge2$? Proposition~\ref{prop:d8-sharp} covers $q<3/2$.
\item In dimensions $d\ge8$, does $p_{2n}$ stay between positive random multiples of its fixed-confidence scale eventually almost surely, or are there fluctuations along subsequences? At $d=8$ the tilt of Section~\ref{sec:heat6} leaves the logarithmic contribution from short time lags in the energy bound: for two endpoints inside the tilted block, the limiting integrand behaves as $h^{-1}$ near zero.
\end{enumerate}

Beyond $\Z^d$, the cut identity of Lemma~\ref{prop:finite-cut} applies on any locally finite graph whose wired forest is a single one-ended tree. The lattice enters at several later steps: the packing argument uses covers of box boundaries, the loop count uses translation invariance and heat-kernel decay, and the contact-count argument interchanges two oriented edges by a lattice reflection. The fluctuation argument also uses independent increments to tilt a segment in a fixed coordinate direction.

\end{document}